\documentclass[a4paper]{amsart}

\usepackage{a4wide}

\usepackage{amsmath}
\usepackage{mathtools}

\usepackage[latin1]{inputenc}
\usepackage[T1]{fontenc}
\usepackage{varioref}
\usepackage{textcomp}
\usepackage{amsfonts}
\usepackage{amssymb}
\usepackage{amsxtra}
\usepackage{stmaryrd} 
\usepackage{bbm}

\usepackage[american]{babel}

\usepackage[all]{xy}
\usepackage{graphicx}

\theoremstyle{definition}
\newtheorem{defin}{Definition}[section]

\newtheorem*{Ack}{Acknowledgements}

\theoremstyle{plain} 
\begingroup 
\newtheorem{thm}[defin]{Theorem}
\newtheorem{prop}[defin]{Proposition}
\newtheorem{lem}[defin]{Lemma}

\newtheorem*{claim}{Claim}
\endgroup

\theoremstyle{remark}
\newtheorem{remark}[defin]{Remark}
\newtheorem*{rem}{Remark}

\newenvironment*{pf}{\textit{Proof}}{\hfill $\Box$ \par}
\newenvironment*{intpf}{\textit{Proof}}{\hfill $\boxtimes$ \par}

\newcommand{\affA}{\mathbbm{A}}
\newcommand{\kxC}{\mathbbm{C}}  
\newcommand{\natN}{\mathbbm{N}} 
\newcommand{\prP}{\mathbbm{P}}  

\newcommand{\F}{\mathcal{F}} 
\newcommand{\G}{\mathcal{G}} 
\newcommand{\II}{\mathcal{I}}
\renewcommand{\IJ}{\mathcal{J}}
\newcommand{\ssh}{\mathcal{O}} 

\newcommand{\liesl}{\mathfrak{sl}}

\newcommand{\lieso}{\mathfrak{so}}

\newcommand{\dZ}{\dot Z}

\newcommand{\Hom}{\operatorname{Hom}}
\newcommand{\id}{\operatorname{id}}
\newcommand{\im}{\operatorname{im}}

\newcommand{\Irr}{\operatorname{Irr}}

\newcommand{\Mat}{\operatorname{Mat}}
\newcommand{\Mor}{\operatorname{Mor}}

\newcommand{\Pf}{\operatorname{Pf}}

\newcommand{\rk}{\operatorname{rk}}

\newcommand{\Spec}{\operatorname{Spec}}

\newcommand{\tensor}{\otimes_{\kxC}}

\newcommand{\red}{\mathbin{\textit{\hspace{-0.8ex}/\hspace{-0.5ex}/\hspace{-0.4ex}}}}
\newcommand{\sred}{\mathbin{\textit{\hspace{-0.8ex}/\hspace{-0.5ex}/\hspace{-0.5ex}/\hspace{-0.4ex}}}}

\newcommand{\mapgets}{\mathrel{\protect\raisebox{0.3ex}{\protect\rotatebox[origin=c]{180}{$\mapsto$}}}}

\newcommand{\acts}{\mathrel{\protect\raisebox{0.3ex}{\protect\rotatebox[origin=c]{270}{$\circlearrowright$}}}}

\newcommand{\GHilb}[1]{G\operatorname{-Hilb}(#1)} 
\newcommand{\dGHilb}[1]{G\operatorname{-Hilb}(#1)^{orb}} 
\newcommand{\SlHilb}[1]{Sl_2\operatorname{-Hilb}(#1)}
\newcommand{\dSlHilb}[1]{Sl_2\operatorname{-Hilb}(#1)^{orb}}
\newcommand{\pHilb}[2]{\operatorname{Hilb}^{#1}(#2)} 

\newcommand{\invHilb}[3]{\operatorname{Hilb}^{#1}_{#2}(#3)} 
\newcommand{\dinvHilb}[3]{\operatorname{Hilb}^{#1}_{#2}(#3)^{orb}}

\newcommand{\finvHilb}[1]{\mathcal{H}ilb_{#1}^G}
\newcommand{\Univ}[3]{\operatorname{Univ}^{#1}_{#2}(#3)} 
\newcommand{\Grass}[2]{\operatorname{Grass}({#1},{#2})} 
\newcommand{\qGrass}[2]{\operatorname{Grass}({#1},{#2})} 
\newcommand{\Grassiso}[2]{\operatorname{Grass}_{iso}({#1},{#2})} 

\allowdisplaybreaks

\author{Tanja Becker} \thanks{This work has been partially supported by DAAD and by the SFB/TR 45 ''Periods, Moduli Spaces and Arithmetic of Algebraic Varieties'' of the DFG}

\title{An example of an $Sl_2$-Hilbert scheme with multiplicities}
\date{October 18, 2010}
\address{Fachbereich Physik, Mathematik und Informatik\\ Johannes Gutenberg -- Universit\"at Mainz \\ D -- 55099 Mainz}

\email{tanja@mathematik.uni-mainz.de}

\begin{document}

\begin{abstract}
We determine the invariant Hilbert scheme of the zero fibre of the moment map of an action of $Sl_2$ on $(\kxC^2)^{\oplus 6}$ as one of the first examples of invariant Hilbert schemes with multiplicities. While doing this, we present a general procedure how to realise the calculation of invariant Hilbert schemes, which have been introduced by Alexeev and Brion in \cite{AB:2005}. We also consider questions of smoothness and connectedness and thereby show that our Hilbert scheme gives a resolution of singularities of the symplectic reduction of the action.
\end{abstract}

\maketitle

\section{Introduction}

Let $G$ be a complex connected reductive algebraic group and $X$ an affine $G$-scheme over $\kxC$. Denote by $\Irr(G)$ the set of isomorphism classes of irreducible representations of $G$ and let $h\colon \Irr(G) \to \natN_0$ be a map, called \textit{Hilbert function} in the following. In this setting, Alexeev and Brion define in \cite{AB:2005} the invariant Hilbert scheme $\invHilb{G}{h}{X}$ parameterising $G$-invariant subschemes of $X$ whose modules of global sections all have the same isotypic decomposition $\bigoplus_{\rho \in \Irr(G)}\kxC^{h(\rho)} \tensor V(\rho)$ as $G$-modules, where $V(\rho)$ denotes the $G$-module corresponding to the irreducible representation $\rho$. This generalises the $G$-Hilbert scheme of Ito and Nakamura \cite{IN:1996}.

In the case where the Hilbert function $h$ is multiplicity-free, i.e. $h(\rho) \in \{0,1\}$ for all $\rho \in \Irr(G)$, several examples of the invariant Hilbert scheme have been determined by Jansou \cite{Jan:2007}, Bravi and Cupit-Foutou \cite{BCF:2008} and Papadakis and van Steirteghem \cite{PvS:2010}, 
which all turn out to be affine spaces. Jansou and Ressayre \cite{JR:2009} give some examples of invariant Hilbert schemes with multiplicities, which are also affine spaces. There are some more involved examples of invariant Hilbert schemes by Brion (unpublished) and Budmiger \cite{Bud:2010}. In this paper, we present a more complex example, namely of an $Sl_2$-Hilbert scheme with Hilbert function 
\begin{equation} \label{hilbfkt} h\colon \natN_0 \to \natN, \;d \mapsto d+1.\end{equation}

The knowledge of such examples where the Hilbert scheme is not an affine space is important for understanding general properties of invariant Hilbert schemes: Which conditions have to be fulfilled so that the invariant Hilbert scheme is connected or smooth? Is the invariant Hilbert scheme a resolution of singularities of the quotient $X \red G$ as the $G$-Hilbert scheme is for finite $G$ up to dimension $3$ \cite{BKR:2001}?

Our example of an $Sl_2$-Hilbert scheme will be smooth and connected and it even will be a resolution of singularities, but it does not inherit the additional structure of symplectic variety of the quotient.

\vspace{1em}

Now let us present the setting of our example. 
Consider the action of $Sl_2$ on $(\kxC^{2})^{\oplus 6} = \Mat_{2 \times 6}(\kxC)$ arising as symplectic double from the action of $Sl_2$ 
on $(\kxC^{2})^{\oplus 3}$ via multiplication on the left.

Let $J = \left(\begin{smallmatrix} 0 & 1 \\ -1 & 0 \end{smallmatrix}\right)$ and $Q = \left(\begin{smallmatrix} 0 & I_3 \\ I_3 & 0 \end{smallmatrix}\right)$.
The moment map $\mu\colon (\kxC^{2})^{\oplus 6} \to \liesl_2$, $M \mapsto MQM^tJ$ defines the symplectic reduction $(\kxC^{2})^{\oplus 6}\sred Sl_2 := \mu^{-1}(0)\red Sl_2$. In \cite{Bec:2010} we obtained its description as a nilpotent orbit closure $\mu^{-1}(0)\red Sl_2 = \overline{\ssh}_{[2^2,1^2]}$ in $\lieso_6$. 
Writing $(\kxC^{2})^{\oplus 6} = \kxC^2 \tensor \kxC^6$ we see that we have a symmetric situation with an action of $SO_6 = SO(Q)$ by multiplication from the right and $\mu$ is invariant for this action, so that $SO_6$ acts on the zero fibre $\mu^{-1}(0)$ and as both actions commute, $SO_6$ also acts on the quotient by $Sl_2$. 
The quotient map $\nu\colon \mu^{-1}(0) \to \mu^{-1}(0)\red Sl_2$ 
is given by mapping $M$ to $M^tJMQ$. In fact, the quotient map of the $Sl_2$-action is the moment map of the $SO_6$-action and vice versa. 
The $SO_6$-action will play an important role while analysing $\mu^{-1}(0)\red Sl_2$ and the corresponding Hilbert scheme.

\vspace{1em}

There are two well-known symplectic resolutions of singularities of the symplectic variety $\overline{\ssh}_{[2^2,1^2]}$, namely the cotangent bundle $T^*\prP^3 \cong \{(A,L) \in Y \times \prP^3 \mid \im A \subset L\}$ and its dual $(T^*\prP^3)^*\cong \{(A,H) \in Y \times (\prP^3)^* \mid H \subset \ker A \}$, where $Y = \{A \in \liesl_4 \mid \rk B \leq 1 \} \cong \overline{\ssh}_{[2^2,1^2]}$. We want to know if there is a distinguished (symplectic) resolution. Since Hilbert schemes of points and $G$-Hilbert schemes are often candidates for (symplectic) resolutions \cite{Fog:1968, Bea:1983, BKR:2001}, we hope that this is also true for invariant Hilbert schemes. Indeed, with the choice of the Hilbert function \eqref{hilbfkt}, in our example we find

\begin{thm} \label{sl2hilbmu}
The invariant Hilbert scheme $\SlHilb{\mu^{-1}(0)} := \invHilb{Sl_2}{h}{\mu^{-1}(0)}$ of the zero fibre of the moment map of the action of $Sl_2$ on $(\kxC^2)^{\oplus 6}$ is the scheme
\begin{equation}\label{result} 
\{(A,W) \in \overline{\ssh}_{[2^2,1^2]} \times \Grassiso{2}{\kxC^6} \mid \im A^t \subset W \}.\end{equation}
It is smooth and connected, and thus a resolution of singularities of the symplectic reduction $\mu^{-1}(0) \red Sl_2$.
\end{thm}

\begin{rem}
$\SlHilb{\mu^{-1}(0)}$ is itself not a \textit{symplectic} resolution of $\mu^{-1}(0) \red Sl_2$, but as it is isomorphic to $\{(A,L,H) \in Y \times \prP^3 \times (\prP^3)^* \mid \im A \subset L \subset H \subset \ker A \}$ via $L = \{v \in \kxC^4 \mid \dim(v \wedge W)=0 \}$, $H = \{v \in \kxC^4 \mid \dim(v \wedge W) \leq 1 \}$ and $W = L\wedge H$, it dominates the two symplectic resolutions:
$$ \begin{xy} \xymatrix{
 & \SlHilb{\mu^{-1}(0)} \ar@{->}[ld] \ar@{->}[dd] \ar@{->}[rd] & \\
T^*\prP^3 \ar@{->}[dr] & & (T^*\prP^3)^* \ar@{->}[dl]\\
 & \mu^{-1}(0) \red Sl_2 &
}\end{xy}
\qquad
\begin{xy} \xymatrix{
 & (A,W) \ar@{|->}[ld] \ar@{|->}[dd] \ar@{|->}[rd] & \\
(A,L) \ar@{|->}[dr] & & (A,H) \ar@{|->}[dl]\\
 & A &}\end{xy} $$
\end{rem}

This paper is organised as follows: In the second chapter we introduce the invariant Hilbert scheme as defined by Alexeev and Brion in \cite{AB:2005}, building upon the work of Haiman and Sturmfels on the multigraded Hilbert scheme \cite{HS:2004}. We give their definition of the invariant Hilbert functor, which is represented by the invariant Hilbert scheme, we introduce the Hilbert-Chow morphism and analyse which conditions on the Hilbert function have to be satisfied so that this morphism, or at least its restriction to a certain component, is proper and birational, the important properties for being a resolution. With regard to this, we define the orbit component $\dinvHilb{G}{h}{X}$, which is the unique component mapping birationally to the set of closed $G$-orbits. If the invariant Hilbert scheme is not irreducible, this component is still a candidate for a resolution.

Afterwards, we turn to our example in chapter $3$. First, we compute the general fibre of the quotient in order to determine the right Hilbert function which guarantees birationality. 

The forth chapter is the heart of this article. First we show how to find generators of the locally free sheaves occurring in the definition of the invariant Hilbert functor in general, then we construct an embedding of the Hilbert scheme into a product of Grassmannians by ideas of Brion and based on the embedding constructed in \cite{HS:2004}. Thus this note not only gives a complex example of an invariant Hilbert scheme with multiplicities of a variety which is not an affine space, but it also can be consulted as a guidance for the determination of further examples.
While describing the general process we always switch to its application to the example at the end of each step. As a result of this, we obtain the orbit component in our example as \eqref{result}.

To conclude the proof of theorem \ref{sl2hilbmu}, i.e. to find out if the orbit component coincides with the whole Hilbert scheme, in chapter $5$ we show that the latter is smooth by considering the tangent space to the invariant Hilbert scheme and we prove that it is connected.

\begin{Ack}
I am very thankful to Michel Brion for introducing me to the world of invariant Hilbert schemes and for guiding me through the determination of this example. I also would like to thank him for his hospitality during the four months I spent in Grenoble. I thank Manfred Lehn and Christoph Sorger for proposing me the work on $G$- and invariant Hilbert schemes and for several discussions about the example. I thank Ronan Terpereau for the exchange of knowledge on invariant Hilbert schemes. I am greatful to Jos\'e Bertin for his private lessons on the $G$-Hilbert scheme which also enlarged my understanding of invariant Hilbert schemes. I greatfully acknowledge the financial support by DAAD and SFB/TR 45.
\end{Ack}

\section{The invariant Hilbert scheme after Alexeev and Brion}

Before passing to the specific example of an invariant Hilbert scheme, we present in general the construction of the invariant Hilbert scheme introduced by Alexeev and Brion in \cite{AB:2004,AB:2005}, which generalizes the $G$-Hilbert scheme for finite groups $G$ after Ito and Nakamura \cite{IN:1996,IN:1999}. For further details on invariant Hilbert schemes consult Brion's survey \cite{Bri:2010}.
\vspace{1ex}

Let $G$ be a complex reductive algebraic group and $X$ an affine $G$-scheme over $\kxC$. Let $\Irr(G)$ denote the set of isomorphism classes of irreducible representations of $G$ and denote by $\rho_0 \in \Irr(G)$ the trivial representation. As $G$ is reductive every $G$-module $W$ decomposes as a sum of its isotypic components
$W = \bigoplus_{\rho \in \Irr G} W_{(\rho)} = \bigoplus_{\rho \in \Irr G} W_{\rho} \tensor V(\rho)$, where $W_{\rho} = \Hom_G(V(\rho), W)$.

We call the dimension of $\Hom_G(V(\rho),W)$ the \textit{multiplicity} of $\rho$ in $W$. If each irreducible representation occurs with finite multiplicity, i.e. for all $\rho \in \Irr G$ we have $h(\rho) := \dim \Hom_G(V(\rho),W) < \infty$, then $h\colon \Irr(G) \to \natN$ is called the \textit{Hilbert function} of $W$.

If $\F$ is a coherent $G$-sheaf over some noetherian basis $S$ where $G$ acts trivially, there is also an isotypic decomposition $\F = \bigoplus_{\rho \in \Irr G} \F_{\rho} \tensor V(\rho)$, where the sheaves of covariants $\F_{\rho} = \mathcal{H}om^G(V(\rho), \F)$ are coherent $\ssh_S$-modules. They are locally free of rank $h(\rho)$ if and only if $\F$ is flat over $S$.

\begin{defin}\cite[Def. 1.5]{AB:2005}
For any function $h\colon \Irr G \to \natN_0$, the associated functor
\begin{align*}
\finvHilb{h}(X)\colon \text{(Schemes)}^{\text{op}} &\to \text{(Sets)}\\
S &\mapsto 
\left\{ 
\begin{array}{c} \xymatrix{
Z \ar@{->}_{p}[dr] & \hspace{-4.2em} \subset &  \hspace{-5.8em} X \times S \ar@{->}^{pr_2}[d]<-2.7em>\\
                   & \hspace{-5em}  &  \hspace{-5.4em} S } \end{array}
\middle| 
\begin{array}{l} Z \text{ a } G\text{--invariant closed subscheme, } \\
p \text{ flat},\\
p_*\ssh_Z \cong \bigoplus_{\rho \in \Irr{G}} \F_{\rho} \tensor V(\rho) \end{array} 
\right\},\\
(f \colon T \to S) & \mapsto (Z \mapsto (id \times f)^*Z)
\end{align*}
such that the sheaves of covariants $\F_{\rho} = \mathcal{H}om^G(V(\rho), p_*\ssh_Z)$ are locally free $\ssh_S$--modules of rank $h(\rho) $, is called the \textit{invariant Hilbert functor}.
\end{defin}

\begin{rem}
In analogy to the case of finite $G$ the coordinate ring of every fibre $Z_s$ of the projection $p\colon Z \to S$ of a closed point $s \in S$ satisfies
$$\kxC[Z_s] = \Gamma(Z_s, \ssh_{Z_s}) = (p_*\ssh_Z)(s) \cong \bigoplus_{\rho \in \Irr{G}} \kxC^{h(\rho)} \tensor V(\rho)$$
since the fibre $\F_{\rho}(s)$ is a $\kxC$-vector space of dimension $h(\rho)$. This can be considered as $h(\rho)$ copies of $V(\rho)$ for every $\rho \in \Irr G$, so we write $\bigoplus_{\rho \in \Irr{G}} h(\rho) V(\rho)$ instead. In particular, the only invariants of $\kxC[Z_s]$ are the elements of the isotypical component of the trivial representation $\rho_0$, i.e. $h(\rho_0)$ copies of the constants.
\end{rem}

\begin{prop}\cite{HS:2004, AB:2004, AB:2005}
There exists a quasi-projective scheme $\invHilb{G}{h}{X}$ representing $\finvHilb{h}(X)$, the {\em invariant Hilbert scheme}.
\end{prop}


There is an analogue of the Hilbert-Chow morphism, the quotient-scheme map
$$\eta\colon \invHilb{G}{h}{X} \to \pHilb{h(\rho_0)}{X\red G}, \; Z \mapsto Z \red G,$$
described in \cite[\S\ 3.4]{Bri:2010}.
It is proper 
and even projective \cite[Prop. 3.12]{Bri:2010}. If we add the condition $h(\rho_0) = 1$, then we have $\eta\colon \invHilb{G}{h}{X} \to \pHilb{1}{X\red G} = X\red G$. We will always assume this in the following.
For birationality one has to choose the Hilbert function $h = h_X$ defined by the isotypic decomposition of the general fibre $F$ of the quotient map $\nu\colon X \to X \red G$:
$$\Gamma(F,\ssh_F) = \bigoplus_{\rho \in \Irr{G}} h_X(\rho)V(\rho).$$

\begin{lem} \label{birat}
If $X$ is irreducible, there is an irreducible component $\dinvHilb{G}{h_X}{X}$ of $\invHilb{G}{h_X}{X}$ such that the restriction of the Hilbert-Chow morphism $\eta\colon \dinvHilb{G}{h_X}{X} \to X\red G$ is birational.
\end{lem}
\begin{pf}
By an independent result of Brion \cite[Prop. 3.15]{Bri:2010} and Budmiger \cite[Thm I.1.1]{Bud:2010}, if $\nu\colon X \to X\red G$ is flat, then $X\red G$ represents the Hilbert functor $\finvHilb{h_X}(X)$, thus $X\red G \cong \invHilb{G}{h_X}{X}$. In the non-flat case let $U \subset X\red G$ be a non-empty open affine subset such that $\nu^{-1}(U) \to U$ is flat. Then $h_{\nu^{-1}(U)} = h_X$ since all fibres of $\nu^{-1}(U) \to U$ have the same Hilbert function as the general fibre of $\nu$, so $U$ is isomorphic to the open subscheme $\invHilb{G}{h_X}{\nu^{-1}(U)} = \eta^{-1}(U)$ of $\invHilb{G}{h_X}{X}$. Thus the restriction of $\eta$ to its closure $\dinvHilb{G}{h_X}{X} := \overline{\eta^{-1}(U)}$ is birational.

If $X$ and hence $X \red G$ is irreducible, so is $U$ and $\eta^{-1}(U) \cong U$. Hence there is an irreducible component $C \subset \invHilb{G}{h_X}{X}$ containing $\eta^{-1}(U)$. The morphism $\eta|_C\colon C \to X\red G$ is dominant and the fibres of an open subset of $X\red G$ are finite (indeed the preimage of each element in $U$ is a point). This means that $\dim C = \dim X \red G$, hence $\overline{\eta^{-1}(U)} = C$ is an irreducible component.
\end{pf}

\begin{defin}
The variety $\dinvHilb{G}{h_X}{X}$ constructed in the lemma is called the \textit{orbit component} or \textit{main component} of $\invHilb{G}{h_X}{X}$. It corresponds to the coherent component for toric Hilbert schemes and is the principal component in the sense that it is birational to the quotient $X \red G$ parameterising the closed orbits of the action of $G$ on $X$.
\end{defin}

\begin{rem} The map $\eta|_{\dinvHilb{G}{h_X}{X}}$ is dominant and proper and $\dinvHilb{G}{h_X}{X} \subset \invHilb{G}{h_X}{X}$ is closed, so $\eta|_{\dinvHilb{G}{h_X}{X}}$ is even surjective.
\end{rem}

\begin{remark}\label{GHilb}
If the general fibre of $\nu\colon X \to X \red G$ happens to be the group $G$ itself, the Hilbert function is $h_X(\rho) = \dim(V(\rho))$ since we have $\Gamma(G,\ssh_G) = \kxC[G] = \bigoplus_{\rho \in \Irr{G}} V(\rho)^* \tensor V(\rho)$ and $\dim(V(\rho)^*) = \dim(V(\rho))$.
In analogy to the case of finite groups we write in this situation
$$\GHilb{X} := \invHilb{G}{h_X}{X} \qquad \text{and} \qquad \dGHilb{X} := \dinvHilb{G}{h_X}{X}.$$
\end{remark}

\section{Determination of the Hilbert function}

\subsection{The quotient related to the Hilbert scheme}

The action of $Sl_2$ on $(\kxC^{2})^{\oplus 3}$ via multiplication on the left is self-dual, so its symplectic double $Sl_2 \times (\kxC^{2})^{\oplus 6} \to (\kxC^{2})^{\oplus 6}$ is also given by multiplication from the left $(g,M) \mapsto gM$. We would like the symplectic structure on $(\kxC^{2})^{\oplus 6}$ to descend to the quotient, so instead of $(\kxC^{2})^{\oplus 6}\red Sl_2$ we consider the symplectic reduction $(\kxC^{2})^{\oplus 6}\sred Sl_2 = \mu^{-1}(0)\red Sl_2$, defined as the quotient of the zero fibre of the moment map $\mu\colon (\kxC^{2})^{\oplus 6} \to \liesl_2$, $M \mapsto MQM^tJ$, where $J = \left(\begin{smallmatrix} 0 & 1 \\ -1 & 0 \end{smallmatrix}\right)$ and $Q = \left(\begin{smallmatrix} 0 & I_3 \\ I_3 & 0 \end{smallmatrix}\right)$. For a more detailed discussion of this action we refer to \cite{Bec:2010}, as well as for the description as a nilpotent orbit closure $\mu^{-1}(0)\red Sl_2 = \overline{\ssh}_{[2^2,1^2]} = \{ A \in \lieso_6 \mid A^2 = 0, \;  \rk A \leq 2,\; \Pf_4(QA) = 0 \}$, where $\Pf_4(QA)$ denotes the Pfaffians of the $15$ skew-symmetric $4 \times 4$-minors of $QA$.
Under the adjoint action this variety consists of two orbits of matrices of rank $2$ and $0$, respectively: $\overline{\ssh}_{[2^2,1^2]} = \ssh_{[2^2,1^2]} \cup \{0\}$.

The quotient map is $\nu\colon \mu^{-1}(0) \to \overline{\ssh}_{[2^2,1^2]}$, $M \to M^tJMQ$.

In coordinates $M = \left(\begin{smallmatrix} x_{11} & x_{12} & x_{13} & x_{14} & x_{15} & x_{16} \\ x_{21} & x_{22} & x_{23} & x_{24} & x_{25} & x_{26} \end{smallmatrix}\right)$ we have
\begin{align*}
M^tJMQ
&= \begin{pmatrix} (-x_{2,i}x_{1,3+j} + x_{1,i}x_{2,3+j} )_{ij} & (-x_{2,i}x_{1,j} + x_{1,i}x_{2,j})_{ij} \\ (-x_{2,3+i}x_{1,3+j} + x_{1,3+i}x_{2,3+j} )_{ij} & (-x_{2,3+i}x_{1,j} + x_{1,3+i}x_{2,j} )_{ij} \end{pmatrix}\\
&= \begin{pmatrix} (\Lambda^{i,3+j} )_{ij} & (\Lambda^{i,j} )_{ij} \\ (\Lambda^{3+i,3+j} )_{ij} & (\Lambda^{j,3+i} )_{ij} \end{pmatrix},
\end{align*}
where $i$ and $j$ always range from $1$ to $3$ and $\Lambda^{s,t} = \det(x^{(s)},x^{(t)})$ is the $2 \times 2$-minor of the $s$-th and $t$-th column in $M$. 
Thus the fibres of $\nu$ consist of those $M$ with fixed $2 \times 2$-minors. A further condition is $M \in \mu^{-1}(0)$, i.e.
\begin{align*}
0 = MQM^t 
&= \begin{pmatrix} 2 \cdot \sum\limits_{i=1}^3 x_{1,i}x_{1,3+i} & \sum\limits_{i=1}^3 (x_{1,i}x_{2,3+i}+x_{1,3+i}x_{2,i}) \\ \sum\limits_{i=1}^3 (x_{1,i}x_{2,3+i}+x_{1,3+i}x_{2,i}) & 2\cdot \sum\limits_{i=1}^3 x_{2,i}x_{2,3+i}\end{pmatrix}.
\end{align*}

\subsection{The general fibre of the quotient} \label{genfib}

In order to determine the Hilbert function $h_{\mu^{-1}(0)}$, so that $\SlHilb{\mu^{-1}(0)} = \invHilb{Sl_2}{h_{\mu^{-1}(0)}}{\mu^{-1}(0)}$ birational to the quotient $\mu^{-1}(0)\red Sl_2$, we have to compute the general fibre of $\nu$. Therefore we need to know the locus where the quotient is flat.
\begin{prop}\label{openorbitflat}
The quotient $\nu$ restricted to the preimage of the open orbit of the $SO_6$-action $\nu^{-1}(\ssh_{[2^2,1^2]}) \to \ssh_{[2^2,1^2]}$ is flat and the fibres over all points in the orbit $\ssh_{[2^2,1^2]}$ are isomorphic.
\end{prop}
\begin{pf}
$\mu^{-1}(0)$ is equipped with an action of $SO_6$ via multiplication on the right, which induces the adjoint action on $\mu^{-1}(0)\red Sl_2$. Since $\nu\colon \mu^{-1}(0) \to \mu^{-1}(0)\red Sl_2 = \overline{\ssh}_{[2^2,1^2]}$ is $SO_6$-equivariant, $\nu$ is flat over the whole $SO_6$-orbit $\ssh_{[2^2,1^2]}$ or over no point of this orbit. By Grothendieck's lemma on generic flatness and $\overline{\ssh}_{[2^2,1^2]}\setminus \ssh_{[2^2,1^2]} = \{0\}$ the second case cannot occur. By equivariance, all fibres over this orbit are  isomorphic.
\end{pf}

\vspace{1em}

As a consequence, for computing the general fibre it is enough to determine the fibre over one point $A_0$ in the flat locus $\ssh_{[2^2,1^2]}$. We choose $A_0 = (a_{ij})$ with $a_{15} = -a_{24} = 1$ and $a_{ij} = 0$ otherwise. For $M \in \nu^{-1}(A_0)$ this corresponds to $\Lambda^{1,2} = 1$, $\Lambda^{i,j} = 0$ otherwise. Thus
$$1 = \Lambda^{1,2} = x_{11}x_{22} - x_{12}x_{21}, \quad \text{hence } x_{11} \neq 0 \neq x_{22} \text{ or } x_{12} \neq 0 \neq x_{21}.$$
Without loss of generality assume $x_{11} \neq 0$. Then $\quad x_{22} = \dfrac{1 + x_{12}x_{21}}{x_{11}}.$

For $j = 3, \hdots, 6$ we have
\begin{align*}
&0 = \Lambda^{1,j} = x_{11}x_{2j} - x_{1j}x_{21} & &\Rightarrow \quad x_{2j} = \frac{x_{1j}x_{21}}{x_{11}},\\
&0 = \Lambda^{2,j} = x_{12}x_{2j} - x_{1j}x_{22} & &\Rightarrow \quad x_{12}\frac{x_{1j}x_{21}}{x_{11}} = x_{1j}\frac{1 + x_{12}x_{21}}{x_{11}} = \frac{x_{1j}}{x_{11}} + \frac{x_{1j}x_{12}x_{21}}{x_{11}}\\
& & &\Rightarrow \quad x_{1j} = 0 \quad \text{for } j = 3, \hdots, 6,\\
& & &\Rightarrow \quad x_{2j} = \frac{x_{1j}x_{21}}{x_{11}} = 0 \quad \text{for } j = 3, \hdots, 6.
\end{align*}
This implies
\vspace{-2.1em} 
\begin{align*}
&x_{11}x_{14} + x_{12}x_{15} + x_{13}x_{16} = 0,\\
&x_{11}x_{24} + x_{12}x_{25} + x_{13}x_{26} + x_{14}x_{21} + x_{15}x_{22} + x_{16}x_{23} = 0,\\
&x_{21}x_{24} + x_{22}x_{25} + x_{23}x_{26} = 0,
\end{align*}
so $M \in \mu^{-1}(0)$ is automatic. This shows that the general fibre is 
$$F := \nu^{-1}(A_0) = 
\left\{ \left(\begin{smallmatrix} x_{11} & x_{12} & 0 & 0 & 0 & 0 \\ x_{21} & x_{22} & 0 & 0 & 0 & 0 \end{smallmatrix}\right) \in (\kxC^2)^{\oplus 6} \; \middle| \; x_{11}x_{22} - x_{12}x_{21} = 1 \right\} \cong Sl_2.$$

\begin{rem}
Analogous calculations over $0$ show that the fibre $\nu^{-1}(0)$ has dimension $5$, so $\nu$ is not flat over $0$ and $\ssh_{[2^2,1^2]}$ is the maximal flat locus.
\end{rem}

\subsection{The Hilbert function of the general fibre}
The Hilbert function is determined by the isotypic decomposition of the general fibre.

The irreducible representations of $Sl_2$ are parametrised by the natural numbers: $\Irr(Sl_2) \cong \natN_0$, $V_d \leftrightarrow d$, where $V_d = \kxC[x,y]_d$ consists of homogeneous polynomials of degree $d$ so that $\dim V_d = d+1$. By remark \ref{GHilb} the coordinate ring of $Sl_2$ decomposes as
$$\kxC[Sl_2] = \bigoplus_{d \in \natN_0} (\dim V_d) V_d = \bigoplus_{d \in \natN_0} (d+1) V_d,$$
so in this case the Hilbert function is given by the dimension $h_{\mu^{-1}(0)}(d) = \dim V_d = d+1$. For the Hilbert scheme this means that the sheaves $\F_d$ have to be locally free of rank $d+1$.

\section{Determination of the orbit component}
Our idea to identify $\SlHilb{\mu^{-1}(0)}$ is to determine generators for the sheaves of covariants $\F_d$ and to use them to embed the $Sl_2$-Hilbert scheme into the product of $\mu^{-1}(0)\red Sl_2$ and some Grassmannian. 
First, in section \ref{shcov} we describe the sheaves $\F_{\rho}$ in general by giving a space of generators $F_{\rho}$ as an $\ssh_{\invHilb{G}{h}{X}}$-module and we calculate $\F_1$ in our example. 
In section \ref{emb} we describe how to obtain a map $\eta_{\rho}$ to the Grassmannian of quotients of $F_{\rho}$ of rank $h(\rho)$ for each $\rho \in \Irr G$. We show that one can embed $\invHilb{G}{h}{X}$ into a product of finitely many of these Grassmannians. Afterwards, for $\SlHilb{\mu^{-1}(0)}$ we calculate the map $\eta_1$ corresponding to the standard representation and we show that this single representation is enough to give an embedding of the orbit component into $\mu^{-1}(0)\red Sl_2 \times \qGrass{F_1}{h(1)}$. Then we determine a strict subset of this which contains the image.
Finally, by writing the Grassmannian as a homogeneous space we prove in section \ref{homogsp} that the embedding is even an isomorphism.
Since the elements of $\SlHilb{\mu^{-1}(0)}$ are subschemes of $\mu^{-1}(0)$, we explicitly determine these subschemes in section \ref{subschemes}.

 \subsection{The sheaves of covariants $\F_{\rho}$} \label{shcov}

To describe the invariant Hilbert scheme or at least its orbit component, we have to determine all possibilities of locally free sheaves $\F_{\rho}$ of rank $h(\rho) $ on $\invHilb{G}{h}{X}$. For the trivial representation we have the following result by Brion \cite[Proof of Prop. 3.15]{Bri:2010}, for which we give a more detailed proof.

\begin{lem}\label{F0} 
If $h(\rho_0) = 1$ then for any scheme $S$ and every subscheme $Z \in \finvHilb{h}(X)(S)$ we have $\F_{\rho_0} = \ssh_S$. In particular, for the universal subscheme $\F_{\rho_0} = \ssh_{\invHilb{G}{h}{X}}$.
\end{lem} 
\begin{pf}
Taking invariants, the defining equation of the $\F_{\rho}$ implies $p_*\ssh_{Z}^G = \bigoplus_{\rho \in \Irr(G)} \F_{\rho} \tensor V(\rho)^G$. But the trivial representation is the only irreducible representation admitting invariants, and all of its elements are invariants. Thus $\bigoplus_{\rho \in \Irr(G)} \F_{\rho} \tensor V(\rho)^G = \F_{\rho_0}$ and there is a morphism $p^{\#}\colon \ssh_{S} = \ssh_{S}^G \to p_* \ssh_{Z}^G = \F_{\rho_0}$ induced by $p$, which is injective since $p$ is surjective.
Both sides are locally free $\ssh_{S}$-modules of rank one.
Over each closed point $s \in S$ the fibres are $\ssh_S(s) = \ssh_{S,s}/\mathfrak{m}_s = k(s) = \kxC$ and $\F_{\rho_0}(s) = (p_*\ssh_Z)^G(s) = (p_*\ssh_Z)^G\tensor k(s) = (p_*\ssh_Z \tensor k(s))^G = \kxC[Z_s]^G$, and $\kxC[Z_s]^G = V(\rho_0) \cong \kxC$. So by Nakayama's lemma, $p^{\#}$ is an isomorphism, hence $\ssh_S \cong \F_{\rho_0}$.
\end{pf}

\vspace{1em}

For general $\rho$, we additionally observe what happens if there is an action on $X$ by another complex connected reductive group $H$ commuting with the $G$-action. By \cite[Prop. 3.10]{Bri:2010}, such an action also induces an action on $X \red G$ and on $\invHilb{G}{h}{X}$, such that the quotient map and the Hilbert-Chow morphism are $H$-equivariant. 

Consider the isotypic decomposition $\kxC[X] = \bigoplus_{\rho \in \Irr G} \kxC[X]_{\rho} \tensor V(\rho)$, where $H$ acts by the induced action on $\kxC[X]_{\rho} = \Hom_G(V(\rho),\kxC[X])$ and trivially on $V(\rho)$.

\begin{prop}\label{genFrho}
For every $\rho \in \Irr G$, the $\kxC[X]^G$-module $\kxC[X]_{\rho}$ is finitely generated, so there is a finite dimensional $H$-module $F_{\rho}$ and an $H$-equivariant surjection $\kxC[X]^G \tensor F_{\rho} \twoheadrightarrow \kxC[X]_{\rho}$. The space $F_{\rho}$ generates $\F_{\rho}$ as an $\ssh_{S}$-module for every scheme $S$ and gives a morphism of $\ssh_S$-$H$-modules $\ssh_S \otimes_{\kxC} F_{\rho} \twoheadrightarrow \F_{\rho}$.
\end{prop}

\begin{pf}
The space $\kxC[X]_{\rho} = \Hom_G(V(\rho),\kxC[X])$ is finitely generated as an $\kxC[X]^G$-module, see \cite[Cor. 5.1]{Dol:2003}. Thus we can choose finitely many generators and define $F_{\rho}$ to be the $H$-module generated by them. This gives an $H$-equivariant surjection $\kxC[X]^G \tensor F_{\rho} \twoheadrightarrow \kxC[X]_{\rho}$.

To determine generators for $\F_{\rho}$ we use the universal subscheme $\Univ{G}{h}{X}$. Then we obtain the result for an arbitrary scheme $S$ and every element in $\finvHilb{h}(X)(S)$ by pulling it back.
We have
$$\xymatrix{
\Univ{G}{h}{X} \ar@{->}_{p}[dr] & \hspace{-4em} \subset &  \hspace{-4em} X \times \invHilb{G}{h}{X} \ar@{->}^{pr_2}[d]<-3em> \ar@{->}[r] & X \ar@{->}^{\nu}[d]\\
                   & \hspace{-4em}  &  \hspace{-4em} \invHilb{G}{h}{X} \ar@{->}^{\eta}[r] & X\red G}$$
The action of $H$ on $X$, $X\red G$ and $\invHilb{G}{h}{X}$ also induces an action of $H$ on $X \times_{X \red G} \invHilb{G}{h}{X}$ and $\Univ{G}{h}{X}$ such that all morphisms in the diagram are $H$-equivariant.
By \cite[Prop. 3.15]{Bri:2010}, the diagram commutes and hence $\Univ{G}{h}{X}$ is even contained in $X \times_{X \red G} \invHilb{G}{h}{X}$.
This inclusion yields a surjective $H$-equivariant morphism $$\ssh_{\invHilb{G}{h}{X}} \otimes_{\kxC[X]^G} \kxC[X] \twoheadrightarrow p_*\ssh_{\Univ{G}{h}{X}}.$$
By definition, we have $p_*\ssh_{\Univ{G}{h}{X}} = \bigoplus_{\rho \in \Irr{G}} \F_{\rho} \otimes_{\kxC} V(\rho)$ with an induced action of $H$ on each $\F_{\rho}$ and the trivial action on $V(\rho)$. Furthermore, we can consider the isotypic decomposition $\ssh_{\invHilb{G}{h}{X}} \otimes_{\kxC[X]^G} \kxC[X] = \bigoplus_{\rho \in \Irr G} \ssh_{\invHilb{G}{h}{X}} \otimes_{\kxC[X]^G} \kxC[X]_{\rho} \tensor V(\rho)$ as $G$-modules.  Together, we obtain $H$-equivariant surjections
$$\ssh_{\invHilb{G}{h}{X}} \otimes_{\kxC[X]^G} \kxC[X]_{\rho} \twoheadrightarrow \F_{\rho}$$
for every $\rho \in \Irr(G)$. 
This shows that the $\ssh_{\invHilb{G}{h}{X}}$-$H$--module $\F_{\rho}$ is generated by $\kxC[X]_{\rho}$, which is in turn generated by $F_{\rho}$ over $\kxC[X]^G$. This yields 
\begin{equation}\label{genmap}
\ssh_{\invHilb{G}{h}{X}} \tensor F_{\rho} 
\twoheadrightarrow \ssh_{\invHilb{G}{h}{X}} \otimes_{\kxC[X]^G} \kxC[X]_{\rho} \twoheadrightarrow \F_{\rho}.
\end{equation}
\end{pf}

\subsection*{Application to $\F_1$}
Now we apply this to our example.
We know that $V_0$ is the trivial representation and by lemma \ref{F0} $\F_0 = \ssh_{\SlHilb{\mu^{-1}(0)}}$ is free of rank $1$.
We suppose that for the representation of lowest dimension $\F_d$ is easiest to compute, so we begin with the standard representation $V_1 = \kxC^2$.
It will turn out in proposition \ref{F1closedemb} that at least the orbit component $\dSlHilb{\mu^{-1}(0)}$ is already completely determined by this sheaf.

There is an action of $SO_6$ on $\mu^{-1}(0)$ via multiplication from the right and the induced action on $\overline{\ssh}_{[2^2,1^2]}$ by conjugation. The induced action on $\SlHilb{\mu^{-1}(0)}$ is also by multiplication from the right. Following proposition \ref{genFrho} we obtain

\begin{prop}
$\F_1$ is generated by the six projections $p_i|_{\mu^{-1}(0)}\colon \mu^{-1}(0) \to \kxC^2$, $i = 1,\dots,6$. Hence we may take $F_1 \cong \kxC^6$ the standard representation of $SO_6$.
\end{prop}

\begin{pf}
Because of proposition \ref{genFrho} and self-duality of the standard representation of $Sl_2$, 
$\F_1$ is generated by $\Hom_{Sl_2}(\kxC^2,\kxC[\mu^{-1}(0)]) = \Mor_{Sl_2}(\mu^{-1}(0),\kxC^2)$. 
The inclusion $\mu^{-1}(0) \subset (\kxC^2)^{\oplus 6}$ induces a surjection 
$\Mor_{Sl_2}((\kxC^2)^{\oplus 6}, \kxC^2) \twoheadrightarrow \Mor_{Sl_2}(\mu^{-1}(0),\kxC^2)$ by shrinking morphisms to $\mu^{-1}(0)$.
By \cite{Ho:1995}, the space of $Sl_2$-equivariant morphisms $\Mor_{Sl_2}((\kxC^2)^{\oplus 6}, \kxC^2)$ is a free module of rank $6$ over the ring of invariants $\kxC[(\kxC^2)^{\oplus 6}]^{Sl_2}$, generated by the projections $p_i\colon (\kxC^2)^{\oplus 6} \to \kxC^2$ to the $i$-th component. 

The restrictions $p_i|_{\mu^{-1}(0)}\colon \mu^{-1}(0) \to \kxC^2$ still span a $6$-dimensional space: Consider for example the matrices $M_i$ where each except the $i$-th column is $0$. Then $M_iQM_i^t = 0$ for $i = 1,\hdots,6$, so $M_i \in \mu^{-1}(0)$. In turn $p_j(M_i) = \delta_{ij}\binom{x_{1j}}{x_{2j}}$ shows that the $p_i|_{\mu^{-1}(0)}$ are linearly independent.
Thus $\Mor_{Sl_2}(\mu^{-1}(0), \kxC^2) \cong \Hom_{Sl_2}((\kxC^2)^{\oplus 6}, \kxC^2)$ and $F_1 = <p_i \mid i = 1,\hdots,6> \cong \kxC^6$.
The $SO_6$-equivariant identification $\kxC^6 \cong \Hom_{Sl_2}((\kxC^2)^{\oplus 6},\kxC^2)$, $e_i \mapsto 
p_i$ induces the inner product $\langle p_i, p_j \rangle = \delta_{i+3,j} + \delta_{j+3,i}$ on $<p_1, \hdots p_6>$. For this reason we can also write $ \langle p, q \rangle = p^tQq$ for all maps $p, q \in F_1$ and we see that $F_1$ is the standard representation.
\end{pf}

 \subsection{Embedding the Hilbert scheme into a product of Grassmannians} \label{emb}

As remarked in the proof of proposition \ref{genFrho}, every map $S \to \invHilb{G}{h}{X}$ gives us a map $\ssh_S \tensor F_{\rho} \to \F_{\rho}$ by pulling back \eqref{genmap}. Since $\F_{\rho}$ is a locally free quotient of $\ssh_S \tensor F_{\rho}$ of rank $h(\rho)$, this in turn corresponds to a map $S \to \qGrass{F_{\rho}}{h(\rho)}$ into the Grassmannian of quotients of $F_{\rho}$ of dimension $h(\rho)$. 
In particular, taking $S = \invHilb{G}{h}{X}$, we obtain a map of schemes
$$\eta_{\rho}\colon \invHilb{G}{h}{X} \to \qGrass{F_{\rho}}{h(\rho)}.$$
In the situation of proposition \ref{genFrho} this map is again $H$-equivariant.
Evaluating at a closed point $s \in S$ yields
\begin{align} \label{corres}
(S &\to \invHilb{G}{h}{X}) & &\longmapsto & &(\ssh_S \tensor F_{\rho} \to \F_{\rho}) & &\longmapsto & &(S \to \qGrass{F_{\rho}}{h(\rho)}),\\
(s &\mapsto Z) & &\longmapsto & &(f_{\rho,s}\colon F_{\rho} \to \F_{\rho}(s)) & &\longmapsto & &(s \mapsto \F_{\rho}(s)),\nonumber
\end{align}
where the fibres $\F_{\rho}(s)$ are vector spaces of dimension $h(\rho)$. Hence for $S = \invHilb{G}{h}{X}$ we have
$$\eta_{\rho}\colon \invHilb{G}{h}{X} \to \qGrass{F_{\rho}}{h(\rho)},\; Z \mapsto \F_{\rho}(Z).$$ 
As $\kxC[X]_{\rho} = \Hom_G(V(\rho),\kxC[X]) \cong \Mor_G(X,V(\rho)^*)$, the elements of the generating space $F_{\rho}$ are $G$-equivariant morphisms from $X$ to $V(\rho)^*$ and evaluating at an element $Z \in \invHilb{G}{h}{X}$ means restricting $\Mor_G(X,V(\rho)^*) \to \Mor_G(Z,V(\rho)^*)$, so in \eqref{corres} we have
$$f_{\rho,Z}\colon F_{\rho} \twoheadrightarrow \F_{\rho}(Z),\; p \mapsto p|_Z.$$

The map $\eta_{\rho_0}$ does not yield any information because $\qGrass{F_{\rho_0}}{h(\rho_0)} = \qGrass{\kxC}{1}$ is only a point.
The product of the Hilbert-Chow morphism and the $\eta_{\rho}$ defines a map
\begin{equation}\label{Grassimmersion}\invHilb{G}{h}{X} \to X\red G \times \hspace{-1em}\prod\limits_{\begin{smallmatrix}\rho \in \Irr(G)\\ \rho \neq 0\end{smallmatrix}}\hspace{-1em} \qGrass{F_{\rho}}{h(\rho)}.\end{equation}

This map is a closed immersion, even if the product ranges over an appropriately chosen finite subset of $\Irr(G)$:
Indeed, let $B = TU$ be a Borel subgroup of $G$, where $T$ is a maximal torus and $U$ the unipotent radical. Assigning to $V(\rho)$ its highest weight gives a one-to-one correspondence between $\Irr G$ and the set of dominant weights $\Lambda^+$ in the weight lattice $\Lambda$ of $T$. Extend $h$ to $\Lambda$ by $0$. Let $V$ be a finite-dimensional $T$-module containing $X \red U$. By \cite[Thm 1.7, Lemma 1.6]{AB:2005}, we have closed embeddings $\invHilb{G}{h}{X} \hookrightarrow \invHilb{T}{h}{X\red U} \hookrightarrow \invHilb{T}{h}{V}$
and each module $\kxC[V]_{\rho}$ is generated by some $\kxC$-vector space $E_{\rho}$ over $\kxC[V]^G$. The $E_{\rho}$ can be chosen as lifts of $F_{\rho}$, so that we have $E_{\rho} \twoheadrightarrow F_{\rho}$ under $\kxC[V] \twoheadrightarrow \kxC[X]$. As shown by \cite[Thm 2.2, 2.3]{HS:2004} the map
$$\invHilb{T}{h}{V} \hookrightarrow \prod_{\rho \in D}\qGrass{E_{\rho}}{h(\rho)}$$
is a closed immersion for an appropriately chosen finite subset $D \subset \Lambda$. Since $h = 0$ outside $\Lambda^+$ we even have $D \subset \Irr(G)$ in our case.
Every quotient of $F_{\rho}$ of dimension $h(\rho)$ is also a quotient of $E_{\rho}$ of dimension $h(\rho)$, so we have an embedding $\qGrass{F_{\rho}}{h(\rho)} \hookrightarrow \qGrass{E_{\rho}}{h(\rho)}$. As every element in $\invHilb{T}{h}{V}$ coming from $\invHilb{T}{h}{X\red U}$ is already generated by $F_{\rho}$, the composite map $\invHilb{T}{h}{X\red U} \hookrightarrow \prod_{\rho \in D}\qGrass{E_{\rho}}{h(\rho)}$ factors through $\qGrass{F_{\rho}}{h(\rho)}$, so that we obtain
$$\begin{xy} \xymatrix{
& \prod\limits_{\rho \in D} \qGrass{F_{\rho}}{h(\rho)} \ar@^{(->}[r]<1mm> & \prod\limits_{\rho \in D} \qGrass{E_{\rho}}{h(\rho)} \\
\invHilb{G}{h}{X} \ar@^{(-->}[ur] \ar@{->>}[d] \ar@^{(->}[r] & \invHilb{T}{h}{X\red U} \ar@^{(-->}[u] \ar@{->>}[d] \ar@^{(->}[r] & \invHilb{T}{h}{V} \ar@^{(->}[u] \ar@{->>}[d]\\
 X \red G \ar@{=}[r] & (X \red U ) \red T \ar@^{(->}[r] & V \red T
}\end{xy}$$

This suggests the following procedure to determine the invariant Hilbert scheme:
Begin with one ``easy'' representation $\rho_i$ and analyse $\eta \times \eta_{\rho_i}$. If this can be shown to be closed immersion, identify the image. Otherwise add another representation and repeat the analysis. This process will stop with some $\eta \times \eta_{\rho_1} \times \hdots \times \eta_{\rho_s}$ being closed immersion.

\subsection*{Determination of $\eta_1$}
The knowledge of $F_1$ gives us an $SO_6$-equivariant map
$$\eta_1\colon \SlHilb{\mu^{-1}(0)} \to \qGrass{F_1}{\dim V_1} = \qGrass{\kxC^6}{2}, \; Z \mapsto \F_1(Z).$$
The fibre $\F_1(Z)$ of the sheaf $\F_1$ is generated by the restrictions of the projections $p_i\colon \mu^{-1}(0) \to \kxC^2$ to the subscheme $Z \subset \mu^{-1}(0)$. 

\begin{prop} \label{F1closedemb}
\begin{enumerate}
 \item The map $\eta \times \eta_1$ is given by
$$\eta \times \eta_1\colon \SlHilb{\mu^{-1}(0)} \to \mu^{-1}(0)\red Sl_2 \times \Grass{2}{\kxC^6}, \; Z \mapsto (Z\red Sl_2, \ker(f_{1,Z})^{\perp}).$$
 \item The image of $\eta \times \eta_1$ restricted to the orbit component $\dSlHilb{\mu^{-1}(0)}$ is contained in $Y := \{ (A,U) \in \overline{\ssh}_{[2^2,1^2]} \times \Grassiso{2}{\kxC^6} \mid \im A^t \subset U\}$.
\end{enumerate}
\end{prop}

\begin{pf}
\textit{1.} To describe the morphism $\eta_1\colon \SlHilb{\mu^{-1}(0)} \to \qGrass{\kxC^6}{2}$ explicitly, we analyse the map $f_{1,Z}\colon F_1 \to \F_1(Z)$. As it is surjective, we have $\F_1(Z) \cong F_1/\ker(f_{1,Z})$. Now we can identify the Grassmannian of quotients with the Grassmannian of subspaces via the canonical isomorphism $\qGrass{\kxC^6}{2} \to \Grass{2}{\kxC^6}$, $F_1/\ker(f_{1,Z}) \mapsto \ker(f_{1,Z})^{\perp}$. Thus $\eta_1$ is the morphism $\eta_1\colon \SlHilb{\mu^{-1}(0)} \to \Grass{2}{\kxC^6}, \; Z \mapsto \ker(f_{1,Z})^{\perp}$.

\textit{2.} Over $\ssh_{[2^2,1^2]}$, we have $\eta \times \eta_1\colon \eta^{-1}(\ssh_{[2^2,1^2]}) \to \ssh_{[2^2,1^2]} \times \Grass{2}{\kxC^6}$, $Z_A \mapsto (A,\ker(f_{1,Z_A})^{\perp}).$
For analysing the image, we choose the special point $A_0 \in \ssh_{[2^2,1^2]}$ again.
The description above shows that $\ker(f_{1,Z_{A_0}})= \hspace{0.1em}<\hspace{-0.1em}p_3, p_4, p_5, p_6>$ with orthogonal complement $\ker(f_{1,Z_{A_0}})^{\perp}= \hspace{0.1em}<\hspace{-0.1em}p_4, p_5>$ by definition of the inner product above.
Since $p_4^tQp_4 = p_4^tQp_5 = p_5^tQp_5 = 0$, this space is isotropic. Thus for every point $A$ in the open orbit, $\ker(f_{1,Z_{A}})^{\perp}$ is isotropic. As being isotropic is a closed condition, $\eta \times \eta_1$ maps the closure of the preimage of $\ssh_{[2^2,1^2]}$ under $\eta$, the orbit component, to the isotropic Grassmannian:
$$\eta \times \eta_1\colon \overline{\eta^{-1}(\ssh_{[2^2,1^2]})} = \dSlHilb{\mu^{-1}(0)} \to \overline{\ssh}_{[2^2,1^2]} \times \Grassiso{2}{\kxC^6}.$$

For the additional condition, we only need to examine $A_0 = \left( \begin{smallmatrix} 
\begin{array}{c} 0\end{array} &
\begin{smallmatrix}0 & 1 & 0 \\ -1 & 0 & 0 \\ 0 & 0 & 0 \end{smallmatrix}\\
\begin{array}{c}0\end{array} & \begin{array}{c}0\end{array}
\end{smallmatrix} \right)$ again.
We can consider 

\vspace{-1em}

$A_0$ and its transpose $A_0^t$ as maps 
\begin{align*}
&A_0\colon F_1 \to F_1, \; p_4 \to -p_2, \; p_5 \to p_1, \; p_i \to 0 \text{ for } i = 1, 2, 3, 6, \\
&A_0^t\colon F_1 \to F_1, \; p_1 \to p_5, \; p_2 \to -p_4, \; p_i \to 0 \text{ for } i = 3, 4, 5, 6.
\end{align*}
Thus we have $\im(A_0^t) = \hspace{0.1em}<\hspace{-0.1em}p_4,  p_5> = \ker(f_{1,Z_{A_0}})^{\perp}$. Since $\eta \times \eta_1$ is $SO_6$-equivariant, the equality $\im(A^t) = \ker(f_{1,Z_{A}})^{\perp}$ holds for every $A$ in the orbit $\ssh_{[2^2,1^2]}$ and we obtain
$$\eta \times \eta_1(\eta^{-1}(\ssh_{[2^2,1^2]}))\subset Y' := \{ (A,U) \in \ssh_{[2^2,1^2]} \times \Grassiso{2}{\kxC^6} \mid \im A^t = U\}.$$
If $A \in \overline{\ssh}_{[2^2,1^2]}\setminus \ssh_{[2^2,1^2]}$, its rank is smaller than $2$ (indeed $A=0$), and so is $\dim (\im A^t)$. So the closure of $Y'$ in $\overline{\ssh}_{[2^2,1^2]} \times \Grassiso{2}{\kxC^6}$ is $Y$.
\end{pf}

We will see in the further examination that $\eta \times \eta_1$ actually is an isomorphism (proposition \ref{SlHilborb}), even on the whole invariant Hilbert scheme (proposition \ref{conn}).

\subsection{The Grassmannian as a homogeneous space} \label{homogsp}
For a further analysis of the image we consider the isotropic Grassmannian as a homogeneous space $\Grassiso{2}{\kxC^6} = SO_6/P$, where $P = (SO_6)_{W_0}$ is the isotropy group of an arbitrary point $W_0 \in \Grassiso{2}{\kxC^6}$. 
We choose $W_0 = \hspace{0.1em}<\hspace{-0.1em}p_1,p_2>$. If $g_{_W} \in SO_6$ is chosen such that $W = g_{_W}W_0$, the isomorphism is
$$\Grassiso{2}{\kxC^6} \to SO_6/P,\; W \mapsto g_{_W} P = [g_{_W}],\; gW_0 \mapgets [g].$$
The existence of the canonical map $f\colon Y \xrightarrow{pr_2} \Grassiso{2}{\kxC^6} \cong SO_6 / P,\; (A, U) \mapsto U \mapsto [g_{_U}]$ shows that $Y$ is an associated $SO_6$-bundle with fibre $E := f^{-1}([I_6]) = pr_2^{-1}(W_0)$:
$$\begin{xy} \xymatrix{
SO_6 \times^P E \ar@{->}[dr] \ar@{->}^{\cong}[rr] &  & Y \ar@{->}^{f}[dl]\\
                   & SO_6/P  & }\end{xy}
$$

\vspace{-6.3em}

$$
\hspace{-4em}
\begin{xy} \xymatrix{
_{(g,A)P} \ar@{|->}[dr] \ar@{|->}[rr] & & _{(gA, gW_0)}\\
                   & ^{\phantom{xm}[g]=gP} & } \end{xy}
$$

\vspace{-6.9em}

$$
\hspace{15em}
\begin{xy} \xymatrix{
&\hspace{-2em} _{\ni\, (A, W)} \ar@{|->}[dl]\\
 _{\ni\, g_{_W}P} & }                                  
\end{xy}
$$

where $SO_6 \times^P E = SO_6 \times E/\hspace{-1ex}\sim$ with $(g,A) \sim (gp^{-1}, pAp^{-1})$.

\begin{lem} The fibre $E = \{ A \in \overline{\ssh}_{[2^2,1^2]} \mid \im A^t \subset W_0 \}$ is one-dimensional.
\end{lem}
\begin{pf} Let $A^t = (a_{ij})$, i.e. $A^tp_i = \sum a_{ji}p_j$. We have
\begin{itemize}
\item $\im A^t \subset W_0 = \hspace{0.1em}<\hspace{-0.1em}p_1,p_2>$, thus $a_{ij} = 0$ if $i = 3,4,5,6$,
\item by duality, $W_0^{\perp} = \hspace{0.1em}<\hspace{-0.1em}p_1,p_2,p_3,p_6> \subset \ker A^t$, which implies $a_{ij} = 0$ if $j = 1,2,3,6$.
\end{itemize}
There only remain $a_{14}$, $a_{24}$, $a_{15}$ and $a_{25}$. But
\begin{itemize}
\item $A^t \in \lieso_6$ implies $a_{14}= a_{25} = 0$ and $a_{24}= -a_{15}$.
\end{itemize}
Thus $E$ is isomorphic to $\affA_{\kxC}^1$.
\end{pf}

\vspace{1ex}

Connecting this to the Hilbert scheme, we have
$$\begin{xy} \xymatrix{
& \mu^{-1}(0)\red Sl_2 & \\
\SlHilb{\mu^{-1}(0)}^{orb} \ar@{->}^{\eta}[ur] \ar@{->}^{f' = f\circ (\eta \times \eta_1)}[dr] \ar@{->}^{\eta \times \eta_1}[rr] &  & Y \cong SO_6 \times^P E \ar@{->}^{pr_1}[ul] \ar@{->}^{f}[dl]\\
                   & SO_6/P  & }                                  
\end{xy}$$
The existence of $f'$ shows, that $\SlHilb{\mu^{-1}(0)}^{orb}$ can be written as an associated $SO_6$-bundle with fibre $F := f'^{-1}([I_6])$ 
and combining the two $SO_6$-bundles we obtain
$$\begin{xy} \xymatrix{
SO_6 \times^P F \ar@{->}[ddr] \ar@{->}_{(\eta \times \eta_1)'}[d] \ar@{->}^{\cong}[rr] &  & \SlHilb{\mu^{-1}(0)}^{orb} \ar@{->}^{\eta \times \eta_1}[d] \ar@{->}_{f'}[ddl]\\
SO_6 \times^P E \ar@{->}[dr] \ar@{->}^{\cong}[rr] &  & Y \ar@{->}^{f}[dl]\\
                   & SO_6/P  & }\end{xy}$$
As $\eta \times \eta_1$ is birational and proper, restricting $(\eta \times \eta_1)'$ to the fibre over a fixed point of $SO_6$ yields a birational and proper morphism $\psi\colon F \to E$. Since $E$ is isomorphic to the affine line, $\psi$ must be an isomorphism. As a consequence:

\begin{prop} \label{SlHilborb} The orbit component of the $Sl_2$-Hilbert scheme is isomorphic to $Y$:
$$\SlHilb{\mu^{-1}(0)}^{orb} \cong \{ (A,U) \in \overline{\ssh}_{[2^2,1^2]} \times \Grassiso{2}{\kxC^6} \mid \im A^t \subset W\}.$$
\end{prop}

\subsection{The points of $\dinvHilb{G}{h}{X}$ as subschemes of $X$} \label{subschemes}

To identify the points of $\dinvHilb{G}{h}{X}$ as subschemes of $X$, we assume there is an embedding
$$\dinvHilb{G}{h}{X} \hookrightarrow X \red G \times \prod_{\rho \in M} \qGrass{F_{\rho}}{h(\rho)},\;
Z \to (Z \red G, (\F_{\rho}(Z))_{\rho \in M})$$
where $M \subset \Irr G$ is a suitable finite subset and $\F_{\rho}(Z) = F_{\rho}/\ker f_{\rho,Z}$ with $f_{\rho,Z}\colon F_{\rho}\to \F_{\rho}(Z)$. This embedding gives us the invariant part and the $\rho$-parts of the ideal $I_Z$ of $Z$ as
\begin{align*}
&(I_Z)^G = I_{Z \red G}\\
&(I_Z)_{\rho} = (\ker f_{\rho,Z}).
\end{align*}
Thus $I_Z \supset I_M := <I_{Z \red G}, \ker f_{\rho,Z} \mid \rho \in M>$. If $I_M$ already has Hilbert function $h$, then $I_Z$ has no further generators and we obtain $I_Z = I_M$.

\subsection*{The points of $\dSlHilb{\mu^{-1}(0)}$ as subschemes of $\mu^{-1}(0)$}

Our next goal is to describe the points $(A,W)$ of $\dSlHilb{\mu^{-1}(0)}$ as subschemes of $\mu^{-1}(0)$. 

\begin{prop}
The subscheme $Z_{A,W} \subset \mu^{-1}(0)$ corresponding to $(A,W) \in Y$ is
$$Z_{A,W} \cong \left\{ \begin{array}{ll} Sl_2, & \text{if } A \in \ssh_{[2^2,1^2]}, \\ \left\{ \left(\begin{smallmatrix} a & b \\ c & d\end{smallmatrix}\right) \middle| ad-bc = 0\right\}, & \text{if } A = 0.\end{array} \right.$$
\end{prop}

\begin{pf}
To show this, we use the embedding
$\eta \times \eta_1\colon \dSlHilb{\mu^{-1}(0)} \to \mu^{-1}(0)\red Sl_2 \times \Grass{2}{\kxC^6}$, $Z \mapsto (Z\red Sl_2,\ker(f_{1,Z})^{\perp})$ 
and we have to compute $Z_{A,W} = (\eta \times \eta_1)^{-1}(A,W)$ or its ideal $I_{A,W}$. 
The action of $SO_6$ on the Hilbert scheme and on $Y$ reduces this to the calculation of one $Z_{A,W}$ for every orbit of $Y$: Since $\eta \times \eta_1$ is $SO_6$-equivariant, all points in the preimage of one orbit are isomorphic. $Y$ decomposes into two $SO_6$-orbits $\{(A,\im A^t)\mid A \in \ssh_{[2^2,1^2]} \} \cong \ssh_{[2^2,1^2]}$ and $\{0\} \times \Grassiso{2}{\kxC^6}$, because the action on $\Grassiso{2}{\kxC^6}$ is transitive.

First we consider $A \in \ssh_{[2^2,1^2]}$. 
Since $\eta$ is an isomorphism of schemes over $\ssh_{[2^2,1^2]}$, we already know that $Z_{A,W} = \eta^{-1}(A) = \nu^{-1}(A) \cong Sl_2$ by section \ref{genfib}.

Now let $A \in \overline{\ssh_{[2^2,1^2]}}\setminus \ssh_{[2^2,1^2]} = \{0\}$. Then $Z_{0,W} \red Sl_2 = 0$, so all $2\times 2$-minors of elements in $Z_{0,W}$ vanish, i.e. $(I_{0,W})^{Sl_2} = (\Lambda_{ij} \mid i,j = 1,\hdots,6)$.
We calculate the subscheme $Z_{0,W}$ explicitly for $W = W_0 := <p_1,p_2>$. Consider $f_{1,Z_{0,W_0}}\colon F_1 \to \F_1(Z_{0,W_0})$, $q \mapsto q|_{Z_{0,W_0}}$. We know that $W_0 = \ker(f_{1,Z_{0,W_0}})^{\perp}$.
If $q = \sum_{i=1}^6 a_i p_i \in \ker(f_{1,Z_{0,W_0}})$, we have $0 = q(M) = \sum_{i=1}^6 a_i \binom{x_{1i}}{x_{2i}}$ for every $M \in Z_{0,W_0}$.
Thus the component of $I_{0,W_0}$ corresponding to the standard representation is $(I_{0,W_0})_1 = (\sum_{i=1}^6 a_i x_{1i},\; \sum_{i=1}^6 a_i x_{2i} \mid q \in W_0^{\perp})$ and for the induces subscheme
$Z'_{0,W_0} := \Spec(\kxC[\mu^{-1}(0)]/((I_{0,W_0})^{Sl_2} + (I_{0,W_0})_1)) \supset Z_{0,W_0}$ we have
$$Z'_{0,W_0} = \{M \in (\kxC^2)^{\oplus 6}\mid MQM^t = 0, \Lambda^{ij} = 0 \; \forall \; i,j, \sum_{i=1}^6 a_i x_{1i} = 0 = \sum_{i=1}^6 a_i x_{2i} \, \forall\; q \in W_0^{\perp}\}.$$
In our case, $W_0^{\perp} = \hspace{0.2em}<\hspace{-0.2em}p_1, p_2, p_3, p_6\hspace{-0.2em}>$, thus letting $q$ be each of these generators yields the equations $x_{1i} = 0 = x_{2i}$ if $i = 1, 2, 3, 6$. This means that $M$ takes the shape $M = \left( \begin{smallmatrix} 0 & 0 & 0 & x_{14} & x_{15} & 0 \\ 0 & 0 & 0 & x_{24} & x_{25}& 0 \end{smallmatrix} \right)$
and $0 = \Lambda^{45} = x_{14}x_{25} - x_{15}x_{24}$. Then the equation $MQM^t = 0$ is automatically fulfilled. So we obtain
$$Z'_{0,W_0} = \left\{ \begin{pmatrix} 0 & 0 & 0 & x_{14} & x_{15} & 0 \\ 0 & 0 & 0 & x_{24} & x_{25}& 0 \end{pmatrix} \in (\kxC^2)^{\oplus 6} \middle| x_{14}x_{25} - x_{15}x_{24} = 0 \right\}.$$
Since this is a flat deformation of $Sl_2$, the corresponding ideal has the correct Hilbert function, which means that $I_{0,W_0} = ((I_{0,W_0})^{Sl_2} + (I_{0,W_0})_1)$ and $Z_{0,W_0} = Z'_{0,W_0}$.
\end{pf}

\section{Properties of the invariant Hilbert scheme}

To determine the whole $Sl_2$-Hilbert scheme, we analyse some of its properties. First, we show that $\SlHilb{\mu^{-1}(0)}$ is smooth in every point of $\dSlHilb{\mu^{-1}(0)}$ in section \ref{smooth}, so that the orbit component is a smooth connected component of $\SlHilb{\mu^{-1}(0)}$. Section \ref{connec} concludes the proof of theorem \ref{sl2hilbmu}, namely $\SlHilb{\mu^{-1}(0)} = \{(A,W) \in \overline{\ssh}_{[2^2,1^2]} \times \Grassiso{2}{\kxC^6} \mid \im A^t \subset W \}$, by showing that $\SlHilb{\mu^{-1}(0)}$ is connected and hence coincides with the orbit component.

\subsection{Smoothness} \label{smooth}

One way to examine smoothness of the Hilbert scheme is to calculate its tangent space at every point. In smooth points the dimension of the tangent space equals the dimension of the Hilbert scheme, in singular points the dimension of the tangent space is bigger. 
In the first case, one concludes that the orbit component is smooth and that there is no additional component of the invariant Hilbert scheme intersecting it, so $\dinvHilb{G}{h}{X}$ is a connected component of the invariant Hilbert scheme.

\vspace{1em}

Let $Z \in \invHilb{G}{h}{X}$, $R := \Gamma(X, \ssh_X)$ and $\II_Z$ the ideal of $Z$ in $\ssh_X$ with space of global sections $I_Z$. 
\begin{prop} \cite[\S\ 1.4]{AB:2005} The tangent space of the Hilbert scheme is given by
$$T_Z\invHilb{G}{h}{X} = Hom^G_R(I_Z,R/I_Z) = Hom^G_{R/I_Z}(I_Z/I_Z^2, R/I_Z) = H^0(\mathcal{H}om^G_{\ssh_Z}(\II_Z/\II_Z^2, \ssh_Z)).$$
\end{prop}

\begin{remark} \label{tangsp} In special situations more can be said about the structure of the tangent space:
\begin{enumerate}
 \item If $Z$ is smooth and contained in the regular part $X_{reg}$ of $X$, then the normal sheaf $\mathcal{N}_{Z/X} := (\II_Z/\II_Z^2)^{\spcheck} = \mathcal{H}om_{\ssh_Z}(\II_Z/\II_Z^2, \ssh_Z)$ is locally free.
This yields a further description of the tangent space
$$T_Z\invHilb{G}{h}{X} = Hom_{R/I_Z}^G(I_Z/I_Z^2, R/I_Z) = H^0(Z,\mathcal{N}_{Z/X})^G.$$
 \item If $Z = Gx \cong G$ is an orbit isomorphic to the group, then we have a commutative diagram
 $$\begin{xy} \xymatrix{
 G \times \mathcal{N}_{Z/X,e} \ar@{->}_{\text{projection}}[dr]<1mm> \ar@{->}^{\sigma}_{\text{action}}[rr] & &\mathcal{N}_{Z/X} \ar@{->}[d]<-3mm> \\
 &   & \hspace{-4em}G \cong Z
 }\end{xy}$$
where $\mathcal{N}_{Z/X,e}$ is the fibre of $\mathcal{N}_{Z/X}$ at the neutral element $e \in G$. The action $\sigma$ restricted to $G \times \mathcal{N}_{Z/X,e}$ is an isomorphism on the fibres: This is clear for $\mathcal{N}_{Z/X,e}$ and true for the other fibres since $\sigma$ is $G$-equivariant. Since both spaces are vector bundles and $\sigma$ is linear, they are isomorphic.

Giving a $G$-invariant section $s\colon G \to G \times \mathcal{N}_{Z/X,e}$ means choosing a point $p \in \mathcal{N}_{Z/X,e}$ such that $s = \id \times p$. This shows
$$T_Z\invHilb{G}{h}{X} = H^0(Z,\mathcal{N}_{Z/X})^G \cong \mathcal{N}_{Z/X,e}.$$
 \item If $Z$ is not smooth we can consider its regular part $Z_{reg}$. If $Z$ is reduced, restricting morphisms to $Z_{reg}$  yields injections
\begin{align*}
& & &Hom_{\ssh_Z}(\II_Z/\II_Z^2, \ssh_Z) \hookrightarrow Hom_{\ssh_{Z_{reg}}}(\II_{Z_{reg}}/\II_{Z_{reg}}^2, \ssh_{Z_{reg}}) \\
&\text{and } & &Hom^G_{\ssh_Z}(\II_Z/\II_Z^2, \ssh_Z) \hookrightarrow Hom^G_{\ssh_{Z_{reg}}}(\II_{Z_{reg}}/\II_{Z_{reg}}^2, \ssh_{Z_{reg}}),
\intertext{
If $Z_{reg} \subset X_{reg}$, taking global sections we obtain}
& & &Hom^G_{R/I_Z}(I_Z/I_Z^2, R/I_Z) \hookrightarrow H^0(Z_{reg}, \mathcal{N}_{Z_{reg}/X_{reg}})^G.
\end{align*}
All these maps are isomorphisms if $Z$ is normal.

\end{enumerate}
\end{remark}

\subsection*{The tangent space of $\SlHilb{\mu^{-1}(0)}$}

In order to find out if the orbit component coincides with the whole Hilbert scheme, we calculate the tangent space to $\SlHilb{\mu^{-1}(0)}$ in every point of $\dSlHilb{\mu^{-1}(0)}$. 

\begin{prop} \label{tangdim}
For every point $Z \in \dSlHilb{\mu^{-1}(0)}$ the dimension of the tangent space is
$$ \dim T_Z\SlHilb{\mu^{-1}(0)} = 6 = \dim \dSlHilb{\mu^{-1}(0)}.$$
Therefore the orbit component is a smooth connected component of the invariant Hilbert scheme.
\end{prop}

\begin{pf}
As before, we only have to consider one point of each $SO_6$-orbit because the dimension of the tangent space is stable in every orbit of the $SO_6$-action. Over the open orbit there is nothing to show, because we know that $\eta^{-1}(\ssh_{[2^2,1^2]}) \cong \ssh_{[2^2,1^2]}$ is smooth.
Over the origin we consider
\begin{align*}
Z := Z_{0,W_0} &= \left\{ \begin{pmatrix} 0 & 0 & 0 & x_{14} & x_{15} & 0 \\ 0 & 0 & 0 & x_{24} & x_{25} & 0 \end{pmatrix} \middle| x_{14}x_{25}-x_{15}x_{24} = 0 \right\} \\
&\cong \left\{ \begin{pmatrix}  \lambda x & \lambda y \\ \mu x & \mu y \end{pmatrix} \middle| x, y \in \kxC, [\lambda:\mu] \in \prP^1 \right\}.
\end{align*}
$Z$ is normal since it is a complete intersection and the codimension of $Z\setminus \dot Z =\{0\}$ in $Z$ is greater than $2$, namely $3$.
We have $Z \subset \mu^{-1}(0)_{sing}$: If $M \in Z$, all of its $2 \times 2$-minors vanish, thus $M \in V(X^tJX) = \mu^{-1}(0)_{sing}$, where $X = \left(\begin{smallmatrix} x_{11} & x_{12} & x_{13} & x_{14} & x_{15} & x_{16} \\ x_{21} & x_{22} & x_{23} & x_{24} & x_{25} & x_{26} \end{smallmatrix}\right)$ describes the coordinates in $\kxC[x_{11},\hdots,x_{26}]$.

From now on we also write $a := x_{14}$, $b := x_{15}$, $c := x_{24}$ and $d := x_{25}$.
Let $\II$ be the ideal of $Z$ in $R := \kxC[\mu^{-1}(0)] = \kxC[x_{11},\hdots,x_{26}]/(XQX^t)$. We have
\begin{align*}
\II &= (x_{11},x_{12},x_{13},x_{16},x_{21},x_{22},x_{23},x_{26},\underbrace{x_{14}x_{25}-x_{15}x_{24}}_{=:z}),\\
R/\II &= \kxC[x_{11},\hdots,x_{26}]/(x_{11},x_{12},x_{13},x_{16},x_{21},x_{22},x_{23},x_{26},x_{14}x_{25}-x_{15}x_{24},XQX^t)\\
&= \kxC[a,b,c,d]/(ad-bc).
\end{align*}

Then $\II/\II^2 = R<x_{11},x_{12},x_{13},x_{16},x_{21},x_{22},x_{23},x_{26},z>$ with relations $XQX^t = 0$:
\begin{align*}
&0 = x_{11}x_{14} + x_{12}x_{15} + x_{13}x_{16} \equiv x_{11}a + x_{12}b\\
&0= x_{11}x_{24} + x_{12}x_{25} + x_{13}x_{26} + x_{14}x_{21} + x_{15}x_{22} + x_{16}x_{23} \equiv x_{11}c + x_{12}d + x_{21}a + x_{22}b\\
&0 = x_{21}x_{24} + x_{22}x_{25} + x_{23}x_{26} \equiv x_{21}c + x_{22}d \hspace{18em} mod\, \II^2.
\end{align*}

\textbf{Reduction to $Z_{reg}$}

We analyse the tangent space $T_Z \SlHilb{\mu^{-1}(0)}$ of the invariant Hilbert scheme by reducing to 
$$\dot Z := Z_{reg} = Z \setminus \{0\} = \{(\lambda v, \mu v) \mid v \in \kxC^2\setminus \{0\}, [\lambda: \mu] \in \prP^1\}.$$
Let $\dot \II$ be the ideal sheaf of $\dot Z$. Then in the situation of remark \ref{tangsp}(3) we even have equality because $Z$ is normal. 
Thus
$$\dim T_Z \SlHilb{\mu^{-1}(0)} = H^0(Z,\mathcal{H}om^G_{\ssh_Z}(\II/\II^2, \ssh_Z))^{Sl_2} = \dim H^0(\dot Z, \mathcal{H}om^G_{\ssh_{\dot Z}}(\dot \II/\dot\II^2, \ssh_{\dot Z}))^{Sl_2}.$$

The open subscheme $\dZ$ of $Z$ is not affine. Since $\dZ = Z\setminus \{0\} = V(ad-bc)\setminus V(a,b,c,d)$, it is covered by the open affine sets $\dZ_a = \Spec R_a$, $\dZ_b = \Spec R_b$, $\dZ_c = \Spec R_c$ and $\dZ_d = \Spec R_d$, where
\begin{align*}
R_a &= (\kxC[a,b,c,d]/(ad-bc))_a = \kxC[a,a^{-1},b,c,d]/(ad-bc) = \kxC[a,a^{-1},b,c],\\
&\quad\text{since $a$ is invertible and $d = \frac{bc}{a}$},\\
R_b &= (\kxC[a,b,c,d]/(ad-bc))_b = \kxC[a,b,b^{-1},d],\\
R_c &= (\kxC[a,b,c,d]/(ad-bc))_c = \kxC[a,c,c^{-1},d],\\
R_d &= (\kxC[a,b,c,d]/(ad-bc))_d = \kxC[b,c,d,d^{-1}].
\end{align*}

To describe the ideal sheaf $\dot \II/\dot \II^2$, we compute it on each set of this covering. As $\dot \II = \II|_{\dZ}$ and $\II$ conincide on an open subset, $\dot \II/\dot \II^2$ is generated by $x_{11},x_{12},x_{13},x_{16},x_{21},x_{22},x_{23},x_{26},z$ with relations
\vspace{-0.8em}
$$\begin{array}{l}
0 = x_{11}a + x_{12}b\\
0 = x_{11}c + x_{12}d + x_{21}a + x_{22}b\\
0 = x_{21}c + x_{22}d.
\end{array}$$

Since $a$ is invertible in $R_a$, the first relation yields $x_{11} = -\frac{b}{a}x_{12}.$
The second relation becomes $0 = -\frac{b}{a}x_{12}c + x_{12}\frac{bc}{a} + x_{21}a + x_{22}b =  x_{21}a + x_{22}b,$ thus $x_{21} = -\frac{b}{a}x_{22}.$ 
Then the third equation $0 = -\frac{b}{a}{x_{22}}c + x_{22}\frac{bc}{a}$ is automatically fulfilled and gives no more information. Denoting $\dot \II_a := \dot \II|_{\dot Z_a}$, this shows that
$$\dot \II_a/\dot \II_a^2 = R_a<x_{12},x_{13},x_{16},x_{22},x_{23},x_{26},z>$$
is free of rank $7$. This means that $\dot \II/\dot \II^2$ is locally free of rank $7$, since we obtain analogously
\begin{align*}
\dot \II_b/\dot \II_b^2 &= R_b<x_{11},x_{13},x_{16},x_{21},x_{23},x_{26},z>,\\
\dot \II_c/\dot \II_c^2 &= R_c<x_{12},x_{13},x_{16},x_{22},x_{23},x_{26},z>,\\
\dot \II_d/\dot \II_d^2 &= R_d<x_{11},x_{13},x_{16},x_{21},x_{23},x_{26},z>.
\end{align*}
If we compute the intersection $\dZ_{ab} = \Spec R_{ab}$ of $\dZ_a$ and $\dZ_b$, we obtain
\begin{align*}R_{ab} &= \kxC[a,a^{-1},b, b^{-1},c,d]/(ad-bc) = \kxC[a,a^{-1},b,b^{-1},c] = \kxC[a,a^{-1},b, b^{-1},d],\\
\dot \II_{ab}/\dot \II_{ab}^2 &= \kxC[a,a^{-1},b,b^{-1},c]<x_{12},x_{13},x_{16},x_{22},x_{23},x_{26},z>\\
&=\kxC[a,a^{-1},b,b^{-1},d]<x_{11},x_{13},x_{16},x_{21},x_{23},x_{26},z>
\end{align*}
with $d = \frac{b}{a}c$ and base change $x_{11} = -\frac{b}{a}x_{12}$ and $x_{21} = -\frac{b}{a}x_{22}.$

\vspace{1em}
\textbf{Reduction of $Sl_2$-linearised sheaves to sheaves linearised w.r.t. a Borel subgroup}

To compute $H^0(\dZ, \mathcal{H}om(\dot \II/\dot \II^2, \ssh_{\dZ}))^{Sl_2}$, we reduce the $Sl_2$-linearised sheaf $\dot \II/\dot \II^2$ on $\dZ$ to a $B$-linearised sheaf on $\kxC^2\setminus \{0\}$, where $B$ is the Borel subgroup of upper triangular matrices of $Sl_2$. 

\begin{claim}
$\dZ$ is an associated $B$-bundle:
$$\dZ = Sl_2 \times^B E, \quad \text{where}\quad
E = \pi^{-1}(e_1) = \{ (\lambda e_1, \mu e_1) \mid [\lambda:\mu] \in \prP^1 \} 
\cong \kxC^2\setminus \{0\}.$$
\end{claim}
\begin{intpf}
There is a natural map
$$\varphi \colon \dZ \to \prP^1 \times \prP^1,\; (\lambda v, \mu v) \mapsto ([v], [\lambda : \mu]).$$
Since $g \cdot (\lambda v, \mu v) = (\lambda gv, \mu gv)$ for every $g \in Sl_2$, $\varphi$ is equivariant for the action $g \cdot ([v], [\lambda : \mu]) = ([gv], [\lambda : \mu])$ on $\prP^1 \times \prP^1$. This yields an equivariant projection
$$\pi \colon \dZ \to \prP^1,\; (\lambda v, \mu v) \mapsto [v].$$
As the action of $Sl_2$ on $\prP^1$ is transitive with isotropy group $B =  \left\{ \left( \begin{smallmatrix}t & u\\ 0 & t^{-1}\end{smallmatrix} \right) \middle|\; t \in \kxC^*, u\in \kxC \right\}$, 
we can write $\prP^1 = Sl_2 / B$.
If $gB \in Sl_2/B$ and $b = \left(\begin{smallmatrix}t & u\\ 0 & t^{-1}\end{smallmatrix}\right)$, we have
$$gb = \begin{pmatrix}g_{11} & g_{12}\\ g_{21} & g_{22}\end{pmatrix}\begin{pmatrix}t & u\\ 0 & t^{-1}\end{pmatrix} = \begin{pmatrix}tg_{11} & ug_{11}+t^{-1}g_{12}\\ tg_{21} & ug_{21}+t^{-1}g_{22}\end{pmatrix},$$
which shows that in the class of $g$, $g_{11}$ and $g_{21}$ are determined up to a scalar and $g_{12}$ and $g_{22}$ can be modified arbitrarily up to the condition $\det(g) = 1$. Therefore the identification $Sl_2/B \cong \prP^1$ is $gB \mapsto [g_{11}: g_{21}]$.

\vspace{1em}

The action of $B$ on $\kxC^2\setminus \{0\}$ induced by the action of $Sl_2$ on $\dZ$ can be computed as follows:
Let $b = \left(\begin{smallmatrix}t & u\\ 0 & t^{-1}\end{smallmatrix}\right) \in B$. Then $be_1 = \binom{t}{0} = te_1$, thus $b(\lambda e_1, \mu e_1) = (t\lambda e_1, t\mu e_1)$ and this means
$$b \cdot (\lambda, \mu) = (t\lambda, t\mu).$$
Hence the action of $B$ on $\kxC^2\setminus \{0\}$ coincides with the action of $\kxC^*$.

Altogether, we have the following commutative diagram
$$\begin{xy}\hspace{2.5em} \xymatrix{
\hspace{-2em} _{(g,(\lambda,\mu))} \ar@{|->}[r] & _{\left(\lambda\binom{g_{11}}{g_{21}}, \mu\binom{g_{11}}{g_{21}}\right)}
}\end{xy}$$

\vspace{-2.5em}

$$\begin{xy}
\hspace{-6em}\xymatrix{_{(g,(\lambda,\mu))} \ar@{|->}[d]\\^{gB}}\hspace{6em}
\xymatrix{
\hspace{-2em} Sl_2 \times^B \kxC^2\setminus \{0\} \ar@{->}[d]<-5mm>
\ar@{->}^{\cong}[r] & \dZ \ar@{->}[d]\\
\hspace{-2em} Sl_2/B \ar@{->}^{\cong}[r] & \prP^1
} \hspace{8em}
\xymatrix{_{(\lambda v, \mu v)} \ar@{|->}[d]\\^{[v]}}
\end{xy}$$

\vspace{-1.5em}

$$\begin{xy} \xymatrix{
\hspace{3em} ^{gB} \ar@{|->}[r] & ^{[g_{11}:g_{21}]}
}\end{xy}$$
\end{intpf}

Now an $Sl_2$-linearised sheaf $\F$ on $\dZ$ corresponds to a $B$-linearised sheaf $\G$ on $\kxC^2\setminus \{0\}$ as well as their duals correspond to each other. If $j\colon \kxC^2\setminus \{0\} \hookrightarrow \dZ$ denotes the inclusion and $e = I_2\cdot B \in Sl_2/B \cong \prP^1$ we obtain $\G$ as the fibre $\F(e) = j^*\F$. 
In the other direction we have $\F = Sl_2 \times^B \G$.

The invariant global sections of corresponding sheaves coincide:
$$H^0(\dZ,\mathcal{H}om_{\ssh_{\dZ}}(\F,\ssh_{\dZ}))^{Sl_2} = H^0(\kxC^2\setminus \{0\},\mathcal{H}om_{\ssh_{\kxC^2\setminus \{0\}}}(\G,\ssh_{\kxC^2\setminus \{0\}}))^{B}.$$
So we take $\F = \II/\II^2$ and are interested in determining the dual of $j^*\F$:

As $\ssh_{\kxC^2} = \kxC[\lambda, \mu]$ and $\kxC^2\setminus \{0\} = \kxC^2\setminus \{0\}_{\lambda} \cup \kxC^2\setminus \{0\}_{\mu}$ the structure sheaf is given by
\begin{align*}
\ssh_{\kxC^2\setminus \{0\}}(\kxC^2\setminus \{0\}_{\lambda}) &= \kxC[\lambda,\lambda^{-1},\mu],\\
\ssh_{\kxC^2\setminus \{0\}}(\kxC^2\setminus \{0\}_{\mu}) &= \kxC[\lambda,\mu,\mu^{-1}].
\end{align*}

In our case the inclusion is 
$j\colon \kxC^2\setminus \{0\} \to \dZ,\; (\lambda,\mu) \mapsto \left(\begin{smallmatrix}\lambda & \mu \\ 0 & 0\end{smallmatrix}\right),$
so on the level of rings we have $a \mapsto \lambda$, $b \mapsto \mu$, $c \mapsto 0$ and $d \mapsto 0$. This means, that $j^*(\dot \II/ \dot \II^2)$ is given by
\begin{align*}
j^*(\dot \II/ \dot \II^2)(\kxC^2\setminus \{0\}_{\lambda}) &= \kxC[\lambda,\lambda^{-1},\mu]<x_{12},x_{13},x_{16},x_{22},x_{23},x_{26},z>,\\
j^*(\dot \II/ \dot \II^2)(\kxC^2\setminus \{0\}_{\mu}) &= \kxC[\lambda,\mu,\mu^{-1}]<x_{11},x_{13},x_{16},x_{21},x_{23},x_{26},z>,\\
j^*(\dot \II/ \dot \II^2)(\kxC^2\setminus \{0\}_{\lambda\mu}) &= \kxC[\lambda,\lambda^{-1},\mu,\mu^{-1}]<x_{12},x_{13},x_{16},x_{22},x_{23},x_{26},z>\\
&= \kxC[\lambda,\lambda^{-1},\mu,\mu^{-1}]<x_{11},x_{13},x_{16},x_{21},x_{23},x_{26},z>
\end{align*}
with base change $x_{11} = -\frac{\mu}{\lambda}x_{12}$ and $x_{21} = -\frac{\mu}{\lambda}x_{22}$.

To compute the dual $j^*(\dot \II/ \dot \II^2)^{\spcheck} = \mathcal{H}om_{\ssh_{\kxC^2\setminus \{0\}}}(\dot \II/ \dot \II^2,\ssh_{\kxC^2\setminus \{0\}})$, denote by $(y_{ij},w)$ the basis dual to $(x_{ij},z)$, i.e. $y_{ij}(x_{kl}) = \delta_{(ij)(kl)}$, $y_{ij}(z) = 0$, $w(x_{ij}) = 0$, $w(z) = 1$. Then we have
\begin{align*}
j^*(\dot \II/ \dot \II^2)^{\spcheck}(\kxC^2\setminus \{0\}_{\lambda}) &= \kxC[\lambda,\lambda^{-1},\mu]<y_{12},y_{13},y_{16},y_{22},y_{23},y_{26},w>,\\
j^*(\dot \II/ \dot \II^2)^{\spcheck}(\kxC^2\setminus \{0\}_{\mu}) &= \kxC[\lambda,\mu,\mu^{-1}]<y_{11},y_{13},y_{16},y_{21},y_{23},y_{26},w>,\\
j^*(\dot \II/ \dot \II^2)^{\spcheck}(\kxC^2\setminus \{0\}_{\lambda\mu}) &= \kxC[\lambda,\lambda^{-1},\mu,\mu^{-1}]<y_{12},y_{13},y_{16},y_{22},y_{23},y_{26},w>\\
&= \kxC[\lambda,\lambda^{-1},\mu,\mu^{-1}]<y_{11},y_{13},y_{16},y_{21},y_{23},y_{26},w>.
\end{align*}
with base change $y_{11} = -\frac{\lambda}{\mu}y_{12}$ and $y_{21} = -\frac{\lambda}{\mu}y_{22}$.

\vspace{1em}

\textbf{Computation of the global sections}

The global sections $H^0(\kxC^2\setminus \{0\},j^*(\dot \II/ \dot \II^2)^{\spcheck})$ are the kernel of the map
\begin{align*}
\varphi \colon H^0(\kxC^2\setminus \{0\}_{\lambda},j^*(\dot \II/ \dot \II^2)^{\spcheck}) \oplus H^0(\kxC^2\setminus \{0\}_{\mu},j^*(\dot \II/ \dot \II^2)^{\spcheck}) &\to H^0(\kxC^2\setminus \{0\}_{\lambda\mu},j^*(\dot \II/ \dot \II^2)^{\spcheck}),\\
(p,q) &\mapsto p|_{\kxC^2\setminus \{0\}_{\lambda\mu}} - q|_{\kxC^2\setminus \{0\}_{\lambda\mu}}.
\end{align*}
Let
\vspace{-2.2em}
\begin{align*}
p &= p_1y_{12}+p_2y_{13}+p_3y_{16}+p_4y_{22}+p_5y_{23}+p_6y_{26}+p_7w, \quad p_i \in \kxC[\lambda,\lambda^{-1},\mu], \phantom{fuelltext}\\
q &= q_1y_{11}+q_2y_{13}+q_3y_{16}+q_4y_{21}+q_5y_{23}+q_6y_{26}+q_7w, \quad q_i \in \kxC[\lambda,\mu,\mu^{-1}].
\end{align*}
Denote $p_i = \frac{p_i^N}{p_i^D}$ and $q_i = \frac{q_i^N}{q_i^D}$ with $p_i^N, q_i^N \in \kxC[\lambda,\mu]$, $p_i^D \in \kxC[\lambda]$ and $q_i^D \in \kxC[\mu]$, $p_i^N, p_i^D$ relatively prime, as well as $q_i^N, q_i^D$.
In $\kxC[\lambda,\lambda^{-1},\mu,\mu^{-1}]$ we have
$$q = -\frac{\lambda}{\mu}q_1y_{12}+q_2y_{13}+q_3y_{16}-\frac{\lambda}{\mu}q_4y_{22}+q_5y_{23}+q_6y_{26}+q_7w.$$
Thus if $i \in \{2,3,5,6,7\}$, for $p$ and $q$ to be equal in $\kxC[\lambda,\lambda^{-1},\mu,\mu^{-1}]$ we must have $p_i = q_i$, i.e. $p_i^N\cdot q_i^D = p_i^D\cdot q_i^N$. As $p_i^N$ and $p_i^D$ have no common factor, $p_i^D$ must divide $q_i^D$. But $p_i^D$ is a polynomial in $\lambda$ while $q_i^D$ is a polynomial in $\mu$. This forces $p_i^D$ to be constant, w.l.o.g. $p_i^D = 1$. This immediately implies $q_i^D = 1$ since $q_i^N$ and $q_i^D$ are coprime. We obtain $p_i^N = p_i = q_i = q_i^N \in \kxC[\lambda,\mu].$

If $i = 1$ or $4$, we see $p_i = -\frac{\lambda}{\mu}q_i$, thus $\mu p_i = -\lambda q_i$. Thus $p_i^N = \lambda \tilde p_i^N$, $q_i^N = -\mu \tilde p_i^N$ and $p_i^D = 1 = q_i^D$ as before.
This yields
\begin{align*}
\ker \varphi = \{(\lambda &p_1y_{12}+p_2y_{13}+p_3y_{16}+ \lambda p_4y_{22}+p_5y_{23}+p_6y_{26}+p_7w,\\
-\mu &p_1y_{11}+p_2y_{13}+p_3y_{16}-\mu p_4y_{21}+p_5y_{23}+p_6y_{26}+p_7w) \mid p_i \in \kxC[\lambda,\mu] \}\\
= \kxC[&\lambda,\mu]<\lambda y_{12},y_{13},y_{16},\lambda y_{22},y_{23},y_{26},w>.
\end{align*}
Thus $H^0(\kxC^2\setminus \{0\}, \mathcal{H}om_{\ssh_{\kxC^2\setminus \{0\}}}(j^*(\dot\II/\dot\II^2),\ssh_{\kxC^2\setminus \{0\}})) = \kxC[\lambda,\mu]<\lambda y_{12},y_{13},y_{16},\lambda y_{22},y_{23},y_{26},w>$ is a free module of rank $7$.

\vspace{1em}

\textbf{Computation of invariants}

Let us now consider the actions of $Sl_2$ and $B$ on these modules. Let $g = \left( \begin{smallmatrix} g_{11} & g_{12} \\ g_{21} & g_{22} \end{smallmatrix} \right)$. Then we have $g\cdot\left( \begin{smallmatrix} x_{1i} \\ x_{2i} \end{smallmatrix} \right) = \left( \begin{smallmatrix} g_{11}x_{1i}+g_{12}x_{21} \\ g_{21}x_{1i}+g_{22}x_{2i} \end{smallmatrix} \right)$, thus
\begin{align*}
g\cdot x_{1i} &= g_{11}x_{1i}+g_{12}x_{2i}, \hspace{4.4cm} g\cdot x_{2i} = g_{21}x_{1i}+g_{22}x_{2i}, \\
g\cdot a &= g_{11}a+g_{12}c, \hspace{5cm} g\cdot c = g_{21}a+g_{22}c, \\
g\cdot b &= g_{11}b+g_{12}d, \hspace{5cm} g\cdot d = g_{21}b+g_{22}d\\
g \cdot z\phantom{_{ij}} &= g(x_{14}x_{25}-x_{15}x_{24})\\ &=(g_{11}x_{14}+g_{12}x_{24})(g_{21}x_{15}+g_{22}x_{25})-(g_{11}x_{15}+g_{12}x_{25})(g_{21}x_{14}+g_{22}x_{24})\\
&= (g_{11}g_{22}-g_{12}g_{21})(x_{14}x_{25}-x_{15}x_{24}) = z.
\end{align*}

The action on the dual is determined by
\begin{align*}
& \left. \begin{array}{l} g\cdot y_{1i}(x_{1i}) = y_{1i}(g^{-1}x_{1i}) = y_{1i}(g_{22}x_{1i}-g_{12}x_{2i}) = g_{22}\\
g\cdot y_{1i}(x_{2i}) = y_{1i}(g^{-1}x_{2i}) = y_{1i}(-g_{21}x_{1i}+g_{11}x_{2i}) = -g_{21}\end{array} \right\} \quad \Rightarrow \quad g \cdot y_{1i} = g_{22}y_{1i} - g_{21}y_{2i},\\
& \left. \begin{array}{l} g\cdot y_{2i}(x_{1i}) = y_{2i}(g^{-1}x_{1i}) = y_{2i}(g_{22}x_{1i}-g_{12}x_{2i}) = -g_{12}\\
g\cdot y_{2i}(x_{2i}) = y_{2i}(g^{-1}x_{2i}) = y_{2i}(-g_{21}x_{1i}+g_{11}x_{2i}) = g_{11}\end{array} \right\} \quad \Rightarrow \quad g \cdot y_{2i} = -g_{12}y_{1i} + g_{11}y_{2i},\\
& \left. \begin{array}{l}g \cdot w (z) = w(g^{-1}z) = w(z) \quad \Rightarrow \quad g \cdot w = w.\end{array} \right.
\end{align*}



Correspondingly, over $\kxC^2\setminus \{0\}$, the action of $g = \left(\begin{smallmatrix}t & u \\ 0 & t^{-1}\end{smallmatrix}\right)$ is
\begin{align*}
&g \cdot \lambda = t\lambda, & &g \cdot x_{1i} = tx_{1i}+ux_{2i}, & &g \cdot y_{1i} = t^{-1}y_{1i},\\
&g \cdot \mu = t\mu,         & &g\cdot x_{2i} = t^{-1}x_{2i},     & &g \cdot y_{2i} = -uy_{1i} + ty_{2i},\\
&			     & &g \cdot z = z, 			  & &g \cdot w = w.
\end{align*}

Now we have $B = TU$ with torus $T = \left\{ \left( \begin{smallmatrix} t & 0 \\ 0 & t^{-1} \end{smallmatrix} \right)\right\}$ and unipotent radical $U = \left\{ \left( \begin{smallmatrix} 1 & u \\ 0 & 1 \end{smallmatrix} \right)\right\}$. Thus we can compute the $B$-invariants stepwise:
$$\kxC[\lambda,\mu]<\lambda y_{12},y_{13},y_{16},\lambda y_{22},y_{23},y_{26},w>^B = (\kxC[\lambda,\mu]<\lambda y_{12},y_{13},y_{16},\lambda y_{22},y_{23},y_{26},w>^U)^T.$$
Let $\underline{u} = \left( \begin{smallmatrix} 1 & u \\ 0 & 1 \end{smallmatrix} \right)$:
\begin{align*}
& \left. \begin{array}{ll}
\underline{u} \cdot \lambda = \lambda \qquad & \underline{u} \cdot \lambda y_{12} = \lambda y_{12} \quad\\
\underline{u} \cdot \mu = \mu         \qquad & \underline{u} \cdot y_{13} = y_{13} \quad\\
\underline{u} \cdot w = w             \qquad & \underline{u} \cdot y_{16} = y_{16} \quad
\end{array} \right\} \text{ invariants}\\
& \left. \begin{array}{l}
\underline{u} \cdot \lambda y_{22} = \lambda (-uy_{12} + y_{22}) = -u \lambda y_{12} + \lambda y_{22}\\
\underline{u} \cdot y_{23} = -uy_{13} + y_{23}\\
\underline{u} \cdot y_{26} = -uy_{16} + y_{26}
\end{array} \right\} \text{ cannot be combined to form invariants.}
\end{align*}

So we have $\quad \kxC[\lambda,\mu]<\lambda y_{12},y_{13},y_{16},\lambda y_{22},y_{23},y_{26},w>^U = \kxC[\lambda,\mu]<\lambda y_{12},y_{13},y_{16},w>.$

\vspace{1em}

To compute the $T$-invariants, let $\underline{t} = \left( \begin{smallmatrix} t & 0 \\ 0 & t^{-1} \end{smallmatrix} \right)$. We obtain
$$\begin{array}{lll}
\begin{array}{l}
\text{degree $1$:} \\
\underline{t} \cdot \lambda = t\lambda\\
\underline{t} \cdot \mu = t\mu\\
\end{array} \qquad \qquad \quad
& \begin{array}{l}
\text{invariants:}\\
\underline{t} \cdot w = w\\
\underline{t} \cdot \lambda y_{12} = t\lambda t^{-1}y_{12} = \lambda y_{12}\\
\end{array} \quad
& \begin{array}{l} 
\text{degree $-1$:}\\
\underline{t} \cdot y_{13} = t^{-1}y_{13}\\
\underline{t} \cdot y_{16} = t^{-1}y_{16}
\end{array}
\end{array}$$
This yields the invariants $w$, $\lambda y_{12}$, $\lambda y_{13}$, $\mu y_{13}$, $\lambda y_{16}$ and $\mu y_{16}$. So we have computed
\begin{align*} 
H^0(\dot Z,\Hom(\dot\II/\dot\II^2,\ssh_{\dot Z}))^{Sl_2} &= H^0(\kxC^2\setminus\{0\},\Hom(\II/\II^2,\ssh_{\kxC^2\setminus\{0\}}))^B \\
&= \kxC<\lambda y_{12}, \lambda y_{13}, \mu y_{13}, \lambda y_{16}, \mu y_{16}, w>. 
\end{align*}

This means that $T_Z \SlHilb{\mu^{-1}(0)}$ is $6$-dimensional and therefore the orbit component of the invariant Hilbert scheme is a smooth connected component.
\end{pf}

\subsection{Connectivity} \label{connec}

In general, the invariant Hilbert scheme can be disconnected. To examine connectivity we look at $\kxC^*$-actions:

If there is a $\kxC^*$-action on $X$ which commutes with the $G$-action, it descends to a $\kxC^*$-action on $X \red G$ so that the quotient map $X \to X \red G$ is $\kxC^*$-equivariant. In this case, one way to find out if the invariant Hilbert scheme is connected is to compute the induced $\kxC^*$-action on $\invHilb{G}{h}{X}$ and to determine all fixed points of $\kxC^* \acts X \red G$. The Hilbert-Chow morphism is proper and $\kxC^*$-equivariant, therefore for every fixed point in the image there is at least one fixed point in every connected component of the fibre of its preimage.

\begin{rem}
Let $(X \red G)_*$ denote the flat locus of the quotient map. Since $\eta|_{\eta^{-1}((X \red G)_*)}$ is an isomorphism, every irreducible component of the invariant Hilbert scheme different from $\dinvHilb{G}{h}{X} = \overline{\eta^{-1}((X \red G)_*)}$ only contains points of the fibres over $X \red G\setminus (X \red G)_*$. If one can show that all connected components of these fibres meet the orbit component, and additionally one knows the orbit component to be smooth, then there cannot be any further component. In this case $\invHilb{G}{h}{X} = \dinvHilb{G}{h}{X}$ is connected.
\end{rem}

\subsection*{Connectedness of $\SlHilb{\mu^{-1}(0)}$}

\begin{prop} \label{conn}
The invariant Hilbert scheme $\SlHilb{\mu^{-1}(0)}$ is connected, hence it coincides with its orbit component and we have
$$\SlHilb{\mu^{-1}(0)} = \dSlHilb{\mu^{-1}(0)} = \{ (A,U) \in \ssh_{[2^2,1^2]} \times \Grassiso{2}{\kxC^6} \mid \im A^t \subset U\}.$$
\end{prop}

\begin{pf}
We consider the action of $\kxC^*$ on $\mu^{-1}(0)$ by scalar multiplication 
and the induced action on $\mu^{-1}(0)\red Sl_2 = \overline{\ssh}_{[2^2,1^2]}$. For $t \in \kxC$ and $M \in \mu^{-1}(0)$ we have $(tM)^tJ(tM)Q = t^2(M^tJMQ)$, thus the action on the quotient is multiplication with $t^2$. Further $A \in \overline{\ssh}_{[2^2,1^2]}$ is invariant 
if and only if $A = 0$, so all fixed points of $\SlHilb{\mu^{-1}(0)}$ map to $0$.

The induced action on $\SlHilb{\mu^{-1}(0)}$ maps $Z$ to $tZ$. If $Z$ is an $Sl_2$-invariant subscheme of $\mu^{-1}(0)$, then $tZ$ is also $Sl_2$-invariant because the action of $Sl_2$ commutes with scalar multiplication. 
Secondly, the global sections of $Z$ and $tZ$ and their isotypic decompositions coincide, 
so indeed $tZ \in \SlHilb{\mu^{-1}(0)}$.

The following lemma shows that set of $\kxC^*$-fixed points in $\SlHilb{\mu^{-1}(0)}$ is $\Grassiso{2}{\kxC^6}$, the fibre of $\dSlHilb{\mu^{-1}(0)}$ over zero. Consequently, $\eta^{-1}(0)$ has no further components, and the same is true for $\SlHilb{\mu^{-1}(0)}$. This shows proposition \ref{conn} and concludes the proof of theorem \ref{sl2hilbmu}.
\end{pf}

\begin{lem}
The set of fixed points in $\SlHilb{(\kxC^2)^{\oplus 6}}$ under the $\kxC^*$-action is isomorphic to the Grassmannian $\Grass{2}{\kxC^6}$. The $\kxC^*$-fixed points in $\SlHilb{\mu^{-1}(0)}$ correspond to the points of $\Grassiso{2}{\kxC^6}$.
\end{lem}

\begin{pf}
Let $Z \subset (\kxC^2)^{\oplus 6}$ be a $\kxC^*$-fixed point in $\SlHilb{(\kxC^2)^{\oplus 6}}$ or $\SlHilb{\mu^{-1}(0)}$, equivalently its corresponding ideal $\II$ is homogeneous. Since then the Hilbert-Chow morphism maps $Z$ to $0$, all $2 \times 2$-minors of each element in $Z$ vanish. Hence $\II$ contains all the $15$ minors $\Lambda^{ij}$.


Now let us analyse the homogeneous invariant ideals $\II$ in $R = \kxC[x_{11},\hdots, x_{26}]$, containing all $\Lambda^{ij}$, with isotypic decomposition $R/\II = \bigoplus_{d \in \natN_0}(d+1)V_d$, where $V_d = \kxC[x,y]_d$ denotes the representation of $Sl_2$ of dimension $d+1$. Afterwards we will restrict to ideals containing $XQX^t$, which correspond to ideals in $R/(XQX^t)$, and which are the fixed points of $\SlHilb{\mu^{-1}(0)}$.

The representation $(\kxC^2)^{\oplus 6} = \Hom(\kxC^6,\kxC^2)$ consists of $6$ copies of $V_1$, so $R = \bigoplus_{n \in \natN_0} S^n(6 V_1)$. 
Since $R = \bigoplus_{n \in \natN_0}S^n(\Hom(\kxC^6,\kxC^2)^*)$ is graded and $\II$ is homogeneous, $R/\II$ is still a graded object, so that we can analyse it by degree. The invariance of $\II$ guarantees that $\II_1$ is a subrepresentation of $\Hom(\kxC^6,\kxC^2)^*$, i.e. there is a sub-vectorspace $V \subset \kxC^6$ such that $\II_1 = \Hom(V,\kxC^2)^*$.
The isotypic decompostion of $R/\II$ requires exactly two copies of $V_1$, and they must already come from $R_1/\II_1$, since no such copy can be contribued or killed by generators of higher degree. If the dimension of $V$ were $5$ or $6$ then $R_1/\II_1$ would be too small and if $\dim V \leq 3$ then $R_1/\II_1$ would be too big.
Thus we know that $\dim V = 4$, so that after a coordinate transformation we can write $\II \supset \IJ = (x_3,y_3,x_4,y_4,x_5,y_5,x_6,y_6,x_1y_2 - y_1x_2)$, since the other $2 \times 2$-minors $x_iy_j - y_jx_i$ do not contribute to the generation of the ideal. Then $R/\IJ \cong \kxC[x_1,y_1,x_2,y_2]/(x_1y_2 - y_1x_2)$ is the coordinate ring of a flat deformation of $Sl_2$ an has isotypic decomposition $\bigoplus_{n \in \natN_0}(n+1)V_n$ as desired. Hence we need no further generators and $\II = \IJ$.

So the fixed points in $\SlHilb{(\kxC^2)^{\oplus 6}}$ under the $\kxC^*$-action correspond to the choice of a $4$-dimensional subspace of $\kxC^6$, which is parameterised by the Grassmannian $\Grass{4}{\kxC^6}$, which coincides with $\qGrass{\kxC^6}{2}$ and $\Grass{2}{\kxC^6}$.

For $Z$ to be contained in $\mu^{-1}(0)$ we have to pick only those ideals which contain $XQX^t$, so that we have $MQM^t=0$ for every $M \in Z$.
We interpret $M \in (\kxC^2)^{\oplus 6}$ as a map $\kxC^6 \to \kxC^2$. 
The fact $M \in Z = \Spec(R/\II)$ means that $M$ vanishes on $V$, so we can interpret it as a map  $\kxC^6/V \to \kxC^2$. As the inner product on $(\kxC^2)^{\oplus 6}$ is induced by the inner product on $\kxC^6$, the condition $MQM^t = 0$ for every $M \in Z$ is equivalent to the vanishing of $v^tQv$ for all $v \in \kxC^6/V$. This show that $\II \supset (XQX^t)$ if and only if $\kxC^6/V$  is an isotropic subspace of $\kxC^6$.

\end{pf}

\begin{rem}
$\SlHilb{(\mu^{-1}(0)}$ is a subscheme of the Hilbert scheme $\SlHilb{(\kxC^2)^{\oplus 6}}$.
The calculation of the fixed points suggests that $\SlHilb{(\kxC^2)^{\oplus 6}}$ contains the whole Grassmannian as the fibre over $0$. Indeed, $\SlHilb{(\kxC^2)^6} = \{((\kxC)^2)^{\oplus 6} \times \qGrass{2}{\kxC^6} \mid \im A^t \subset W\}$ as forthcoming work by Terpereau will show.
\end{rem}

\end{document}